\documentclass[a4paper,11pt]{article}
\usepackage{}
\usepackage{bbm}
\usepackage{CJK}
\usepackage{latexsym}
\usepackage{graphicx}
\usepackage{subcaption}
\usepackage{amssymb}
\usepackage{amsmath}
\usepackage{amssymb,amsmath,amsthm}
\usepackage[numbers,sort&compress]{natbib}
\usepackage[colorlinks,linkcolor=blue,anchorcolor=blue,citecolor=blue]{hyperref}
\usepackage[colorinlistoftodos]{todonotes}
\usepackage{tikz}
\usepackage{amsmath}
\usepackage{array}
\usepackage{booktabs}
\usepackage{graphicx}

\newcolumntype{L}[1]{>{\centering\arraybackslash}m{#1}}

\renewcommand{\thefootnote}{\fnsymbol}

\newtheorem{theorem}{Theorem}[section]
\newtheorem{corollary}{Corollary}[section]
\newtheorem{lemma}{Lemma}[section]
\newtheorem{example}{Example}[section]
\newtheorem{definition}{Definition}[section]
\newtheorem{remark}{Remark}[section]
\newtheorem{proposition}{Proposition}[section]
\newtheorem{condition}{Condition}[section]
\usepackage{amsmath, mathtools, extarrows}
\usepackage{booktabs}
\usepackage{enumerate}
\usepackage{kantlipsum, tikz}
\usepackage{amsmath, mathtools, extarrows}
\usepackage{booktabs}

\usetikzlibrary{arrows.meta, positioning, calc, patterns} 

\large\normalsize

\begin{document}

\title{Decomposition Method for Lipschitz Stability of General LASSO-type Problems%
\footnotemark[1]}
\author{Chunhai Hu \footnotemark[2]\qquad  Wei Yao \footnotemark[3]\qquad Jin Zhang \footnotemark[4]}
\renewcommand{\thefootnote}{\fnsymbol{footnote}}
\footnotetext[2]{Department of Mathematics, Southern University of Science and Technology, Shenzhen {\rm518055}, China (huch@sustech.edu.cn, chunhaihu@outlook.com).}
\footnotetext[3]{National Center for Applied Mathematics Shenzhen, Department of Mathematics, Southern University of Science and Technology, Shenzhen {\rm518055}, China (yaow@sustech.edu.cn).}
\footnotetext[4]{Corresponding author. Department of Mathematics, Southern University of Science and Technology, National Center for Applied Mathematics Shenzhen, Shenzhen {\rm518055}, China (zhangj9@sustech.edu.cn).}
\date{}
\maketitle
\begin{abstract}
This paper introduces a decomposition-based method to investigate the Lipschitz stability of solution mappings for general LASSO-type problems with convex data fidelity and $\ell_1$-regularization terms.
The solution mappings are considered as set-valued mappings of the measurement vector and the regularization parameter.
Based on the proposed method, we provide two regularity conditions for Lipschitz stability: the weak and strong conditions.
The weak condition implies the Lipschitz continuity of solution mapping at the point in question, regardless of solution uniqueness.
The strong condition yields the local single-valued and Lipschitz continuity of solution mapping.
When applied to the LASSO and Square Root LASSO (SR-LASSO), the weak condition is new, while the strong condition is equivalent to some sufficient conditions for Lipschitz stability found in the literature.
Specifically, our results reveal that the solution mapping of the LASSO is globally (Hausdorff) Lipschitz continuous without any assumptions.
In contrast, the solution mapping of the SR-LASSO is not always Lipschitz continuous. A sufficient and necessary condition is proposed for its local Lipschitz property.
Furthermore, we fully characterize the local single-valued and Lipschitz continuity of solution mappings for both problems using the strong condition.
\end{abstract}


{\bf AMS}
49J53, 49K40, 62J07, 90C25, 90C31

\section{Introduction}
\label{introduction-1}

In this work, we consider the following LASSO-type problem with general convex data fidelity and $\ell_1$-regularization terms:
\begin{equation}\label{Psc}
\text{\textbf{(The general problem)}}\quad\quad   \min \limits_{x\in \mathbb{R}^n} h(Ax-b)+\lambda \|x\|_1,
\end{equation}
where $A\in \mathbb{R}^{m\times n}$ is a measurement (or sensing) matrix, $b\in \mathbb{R}^m$ is a vector of measurements (or observations), $\lambda>0$ is a regularization (or tuning) parameter, and $h:\mathbb{R}^m\rightarrow \mathbb{R}_+$ is a general loss function assumed to be continuous and convex. 
The LASSO-type problem and its variant have a wide range of applications in 
statistics \cite{Tib96, Belloni14, Lockhart14, Tian18, Tibshirani11, Geer16, Zou06}, compressed sensing \cite{Adcock19, Adcock22b, Foucart23}, signal processing \cite{Adcock19, Adcock21a, Adcock22a, Adcock22b, Foucart23, Hastie09} and machine learning \cite{Adcock21b, Colbrook22, Hastie09}. 
For example, they are introduced in \cite{Tib96,Belloni11} as optimization models to recover a sparse signal $x_0\in\mathbb{R}^n$ from noisy linear measurements $b=A x_0 + w$, where $w\in\mathbb{R}^m$ is a noise (or perturbation) vector. 

The setting of loss functions in \eqref{Psc} covers many important cases, including  $h(z):=\| z\|_p^q/q$ with $p,q\geq1$ (see \cite{Foucart23,binev2024optimal} and the references therein).
Specifically, the LASSO-type problem \eqref{Psc} generalizes the least absolute deviation LASSO problem \cite{wang2007robust} ($p=1,q=1$), the LASSO problem \cite{Tib96} ($p=2,q=2$):
\begin{equation}\label{P1}
\text{\textbf{(LASSO)}}\quad\quad\quad  \min \limits_{x\in \mathbb{R}^n}\frac{1}{2}\|Ax-b\|^2+\lambda \|x\|_1,
\end{equation}
and the Square Root LASSO (SR-LASSO) problem \cite{Belloni11} ($p=2,q=1$):
\begin{equation}\label{Pi1}
\text{\textbf{(SR-LASSO)}}\quad\quad\quad  \min \limits_{x\in \mathbb{R}^n}\|Ax-b\|+\lambda \|x\|_1. 
\end{equation}

For fixed $A$, $b$ and $\lambda$,  
a wide variety of numerical methods have been developed to solve the LASSO-type problem (see, \textit{e.g.}, \cite{beck2009fast, lee2014proximal, kanzow2021globalized, milzarek2014semismooth} and the references therein).
However, in some important situations, we must consider scenarios where both \( b \) and \( \lambda \) vary naturally. 
For example, this occurs when the vector \( b \) obtained contains inaccuracies or when the tuning parameter \( \lambda \) needs optimal selection \cite{BBH,BBH2,nghia2024geometric,meng2024lip}. From the perspective of perturbation analysis, the LASSO-type problem \eqref{Psc} is an optimization problem with parameters. Therefore, it is of great importance to investigate the behavior of solution mapping when both $b$ and $\lambda$ are perturbed. 

Although there are several sufficient conditions to guarantee solution uniqueness for the LASSO and SR-LASSO  \cite{tibshirani2013lasso,ZY,gilbert2017solution,BBH,BBH2,nghia2024geometric}, it is still difficult to avoid cases of multiple solutions. Solution non-uniqueness happens even if we assume that $A$ is of full row rank and the tuning parameter $\lambda$ is optimal. We provide an example to illustrate this. 

\begin{example}\label{fla}
We consider the following sparse recovery via $\ell_1$ minimization in the noisy setting based on the LASSO and SR-LASSO: 
$$A=\begin{pmatrix}
1&3&0\\
1&-1&1
\end{pmatrix}\in \mathbb{R}^{2\times 3}, ~x_0=\begin{pmatrix}
0\\
0\\
1
\end{pmatrix}\in\mathbb{R}^{3},
~b=Ax_0+
\gamma \begin{pmatrix}
0\\
1
\end{pmatrix}=\begin{pmatrix}
0\\
1+\gamma
\end{pmatrix}\in\mathbb{R}^{2},$$
where $\gamma$ is the noise scale. 
For the LASSO problem \eqref{P1}, the optimal solution set is
\begin{equation*}
    \aligned
S_{LA}(b,\lambda)
=&\left\{(-3s,s,1+\gamma-\lambda+4s)^T\,:
\,\frac{\lambda-1-\gamma}{4} \leq s \leq 0 \right\}~~\forall \lambda \in \left(0,1+\gamma\right),
    \endaligned
\end{equation*}
and
$S_{LA}(b,\lambda)=\{0\}\subset \mathbb{R}^3$ for all $\lambda\in \left[1+\gamma,\infty\right)$.  
Then the best parameter choice
$$\lambda^*_{LA}:=\arg\min\limits_{\lambda>0}d(x_0,S_{LA}(b,\lambda))=\gamma \in \left(0,1+\gamma\right) .$$
Note that $S_{LA}(b,\lambda)$ is not single-valued for all $\lambda>0$ around the best parameter $\lambda^*_{LA}$.

For the SR-LASSO problem \eqref{Pi1}, the optimal solution set is
\begin{equation*}
    \aligned
S_{SR}(b,\lambda)
=&\left\{(-3s,s,1+\gamma+4s)^T\,:
\,\frac{-1-\gamma}{4} \leq s \leq 0 \right\}~~\forall \lambda\in \left(0,1\right),
    \endaligned
\end{equation*}
\begin{equation*}
    \aligned
S_{SR}(b,1)
=&\left\{(-3s,s,t)^T\,:
\,\frac{-1-\gamma}{4} \leq s \leq 0,\, 0\leq t\leq 1+\gamma+4s \right\}~~\mbox{for}~\lambda =1,
    \endaligned
\end{equation*}
and
$S_{SR}(b,\lambda)=\{0\}\subset \mathbb{R}^3$ for all $\lambda\in \left(1,\infty\right)$. 
Then the best parameter choice
$$\lambda^*_{SR}:=\arg\min\limits_{\lambda>0}d(x_0,S_{SR}(b,\lambda))=1.$$
Clearly, $S_{SR}(b,\lambda)$ is not single-valued for all $0<\lambda\leq\lambda^*_{SR}$.
  \end{example}
Therefore, it is essential to study the behavior of solution mapping for \eqref{Psc}, regardless of whether the solution is unique, i.e., in the multi-valued case.  
The Lipschitz property of set-valued mappings is particularly important for analyzing the behavior of feasible sets and optimal solution sets in optimization problems when parameters change \cite{Shapiro00,DR09}. 
In fact, it has received much attention in the last decades, we refer the interested reader to \cite{Hoffman52, Mangasarian87,Klatte87,liwu} for parametric optimization in linear or polyhedral cases; \cite{Yen95, Robinson07} for variational inequalities in polyhedral cases; \cite{Lewis09} for robust regularization problems; \cite{mordukhovich2014} for second-order cone programs;  \cite{xu2016, chenjiang2018, guoxu2023} for stochastic programming; and \cite{bolte2024} for parametric monotone inclusion problem. 
For comprehensive references, see the books by Rockafellar and Wets \cite{rockwets98}, Bonnans and Shapiro \cite{Shapiro00}, Mordukhovich \cite{M2006, M2018}, Dontchev and Rockafellar \cite{DR09} and Ioffe \cite{Ioffe2017}.

Let $S(b,\lambda)$ be the solution mapping of \eqref{Psc}, considered as a set-valued mapping of the measurement vector $b$ and the regularization parameter $\lambda$. 
The purpose of this paper is to develop a method for the Lipschitz stability of $S(b,\lambda)$ in both multi-valued and single-valued cases. 
This type of Lipschitz stability was studied in \cite{BBH} for the LASSO and in \cite{BBH2} for the SR-LASSO, focusing on single-valued solution mappings. Specifically, Berk \textit{et al.} \cite{BBH,BBH2} established sufficient conditions for the single-valued Lipschitz and smoothness properties of solution mappings by proving the strong metric regularity of the subdifferentials of  objective functions. 
They also provided explicit formulas for the Lipschitz modulus. Recently,  
Nghia \cite{nghia2024geometric} fully characterized the local single-valued and Lipschitz property of solution mappings for the group LASSO problem and the matrix setting of the LASSO problem \eqref{P1} with the nuclear norm replacing the $\ell_1$-norm. 
Nghia's approach utilizes the recent geometric characterizations of tilt stability from \cite{nghia2024geometric} and strong minima from \cite{fadili2023geometric}. 
Note that all of these works focus on single-valued case.   
To the best of our knowledge, only recently Meng et al. \cite{meng2024lip} obtained the Lipschitz continuity of solution mapping for the LASSO in the multi-valued case, assuming that $A$ has full row rank. They established a verifiable sufficient condition for a polyhedral mapping to be Lipschitz continuous on its domain and then applied this condition to the LASSO and basis pursuit problems.

\textbf{Our approach and novelty.} We approach Lipschitz stability differently than the methods in \cite{BBH,BBH2,nghia2024geometric,meng2024lip}. These works investigated the Lipschitz stability of $S(b,\lambda)$ by directly studying the subdifferentials of objective functions. 
We present an alternative procedure. 
Our approach decomposes the analysis of the Lipschitz stability of $S(b,\lambda)$ into two stages: 
\begin{enumerate}[(I)]
    \item Establish the single-valued and Lipschitz property of the mapping 
    \begin{equation}\label{1521a-0}
 H(b,\lambda):=(A\circ S) (b,\lambda)=\{Ax:x\in S(b,\lambda)\}~~\forall (b,\lambda)\in  \mathbb{R}^m\times \mathbb{R}_{++} ,
\end{equation}
under proper conditions; 
    \item Prove the Lipschitz continuity of $S$ using the established property of $H$.
\end{enumerate}
The motivation for our approach is as follows: (1) if $S$ is Lipschitz continuous, then $H$ is also Lipschitz continuous; (2) it is well known that $\{Ax:x\in S(b,\lambda)\}$ is a singleton for the LASSO. Hence, it is reasonable to expect that $H$ is single-valued and Lipschitz continuous under proper conditions. 
However, each task in the two stages is challenging. 
It’s unclear whether studying the Lipschitz stability of $H$ is easier than that of $S$. 
One important observation in our study is that the mapping $H$ can be represented as 
\begin{equation}\label{H-represent}
    H(b,\lambda)=\left\{y\in \mathbb R^m:0\in \frac{1}{\lambda}\partial h(y-b) +  G(y)\right\}
    ~~\forall (b,\lambda)\in  \mathbb{R}^m\times \mathbb{R}_{++} ,
\end{equation}
where $G:\mathbb{R}^m\rightrightarrows \mathbb{R}^m$ is a maximal monotone mapping with its inverse mapping
\begin{equation}\label{G-inverse}
    G^{-1}(z)= A(\partial \|\cdot\|_1)^{-1}(A^Tz)
    ~~\forall z\in \mathbb R^m.
\end{equation}
Thanks to the representation \eqref{H-represent}, we can establish the single-valued and Lipschitz property of $H$ under proper conditions. For instance, we can do this by using the coderivative criterion for generalized equations in \cite[Exercise 4C.5]{DR09}, which is essentially based on Mordukhovich criterion (see \cite[Theorem 3.3]{M2018}, \cite[Theorem 9.40]{rockwets98}). 

After establishing the property of $H$ in Stage I, proving the result in Stage II remains challenging. 
A satisfactory method to recover $S$ from $H$ that allows Lipschitz continuity to transfer from $H$ to $S$ is unknown, even if $H$ is single-valued. 
We address this issue by decomposing $S$ into the union of finitely many simpler set-valued mappings defined in \eqref{510g}, namely, 
\begin{equation}\label{S-decomposition}
S(b,\lambda)=\bigcup\limits_{J_1,J_2\subset  [n],~J_1\cap J_2=\emptyset}S_{J_1 J_2}(b,\lambda)~~\forall (b,\lambda)\in \mathbb{R}^m\times \mathbb{R}_+ .
\end{equation}
More importantly, when $H$ is single-valued on a set $U$, we can prove that 
\begin{equation}\label{S-J1J2-F}
     S_{J_1 J_2}(b,\lambda)=\left( F_{J_1 J_2}\circ H \right) (b,\lambda)
     ~~\forall (b,\lambda)\in \mathrm{dom}(S_{J_1 J_2})\cap U,
\end{equation}
where $F_{J_1 J_2}:\mathbb R^m\rightrightarrows \mathbb{R}^n$ is an elementary polyhedral multifunction defined by 
$$F_{J_1 J_2}(y):=\{x\in P_{J_1 J_2}:Ax=y\}~~\forall y\in \mathbb R^m.$$
Here $P_{J_1 J_2}$ is a simple polyhedral set defined in \eqref{ePE}. 
Hence, we can obtain the Lipschitz property of $S_{J_1 J_2}$ from the single-valued and Lipschitz property of $H$. 
However, we cannot directly derive the Lipschitz continuity of $S$ from the decomposition \eqref{S-decomposition}. Instead, we can obtain its outer Lipschitz continuity. 
Based on a notable result on Lipschitz continuity from Robinson \cite[Theorem 1.5]{Robinson07}, the next step is to show the inner semicontinuity of $S$. This follows from an important consistency property of the decomposition \eqref{S-decomposition}, as stated in Proposition \ref{lemma-S-lip-inner}.

\textbf{Main contributions.} 
We propose a decomposition-based method to investigate the Lipschitz stability of  solution mapping $S(b,\lambda)$ of \eqref{Psc}, in both multi-valued and single-valued cases.  
Based on the proposed method, our first main contribution concerns the Lipschitz continuity of $S$, regardless of whether the solution is unique. Concretely, given $(\bar b,\bar\lambda)\in \mathbb R^m\times \mathbb R_{++}$, we establish in Theorem \ref{T01} that $S$ is locally (Hausdorff) Lipschitz continuous at $(\bar b, \bar \lambda)$ if the following condition holds:
\begin{condition}(Weak)\label{conditionH1-0}
There is $\bar x\in S(\bar b, \bar \lambda)$ such that $h$ is twice continuously differentiable around $A\bar x-\bar b$ and $\ker(\nabla^2h(A\bar x- \bar b))\cap \mathrm{rge}(A_{J_h})=\{0\}$, where 
\begin{equation}\label{J-LASSO}
    J_{h}=J_{h}(\bar x):=\left\{i\in [n]:|A_i^T\nabla h(A\bar x-\bar b)|=\bar\lambda\right\}.
\end{equation}
\end{condition}
The independence of $J_{h}$ from $\bar x$ arises because $H(b,\lambda)$ is a singleton under Condition \ref{conditionH1-0}, as proven in Proposition \ref{L514a}. 
It is worth noting that Condition \ref{conditionH1-0} holds automatically for strongly convex loss functions because  $\ker(\nabla^2h(A\bar x- \bar b))=\{0\}$.

Our second main contribution addresses the local single-valued and Lipschitz continuity of $S$.  
We show that the following condition is sufficient for these properties.
\begin{condition}(Strong)\label{conditionH2-0}
Condition \ref{conditionH1-0} holds, and  $A_{J_{h}}$ has full column rank.
\end{condition}
Specifically, in Theorem \ref{T318-h}, given $(\bar b,\bar\lambda)\in \mathbb R^m\times \mathbb R_{++}$, we prove that $S$ is locally single-valued and Lipschitz continuous around $(\bar b,\bar \lambda)$  if Condition \ref{conditionH2-0} holds.

Next we apply our theory to the LASSO and SR-LASSO.  
First, by leveraging the polyhedral structure of the LASSO and using a landmark result of Robinson \cite{Robinson07} on polyhedral multifunctions, we prove in Theorem \ref{Ti1} that the solution mapping $S$ of the LASSO is globally (Hausdorff) Lipschitz continuous without any assumptions. This constitutes our third main contribution. 
Furthermore, in Theorem \ref{T318-lasso}, 
given $(\bar b,\bar\lambda)\in \mathbb R^m\times \mathbb R_{++}$, 
we fully characterize the local single-valued and Lipschitz continuity of $S$ around $(\bar b,\bar \lambda)$ using the condition $\mathrm{ker}(A_J)=\{0\}$. 
Here $J:=\{i\in [n]:|A_i^T( \bar b-A\bar x)|= \bar \lambda\}$. This condition, introduced in \cite{tibshirani2013lasso} for solution uniqueness, is known as Tibshirani's sufficient condition. Berk \textit{et al.} \cite{BBH} proved that it is also a sufficient condition for local Lipschitz continuity.
It is worth noting that recently Nghia \cite[Theorem 3.12]{nghia2024geometric} also proved that Tibshirani's sufficient condition is necessary, 
using the second subderivative characterization of strong minima from \cite[Theorem 13.24]{rockwets98}.
In contract, our approach uses only the necessary result of solution uniqueness for the LASSO from Zhang \textit{et al.} \cite[Theorem 2.1]{ZY}. 

Now let's consider the SR-LASSO problem \eqref{Pi1}. 
Unlike the LASSO, the solution mapping $S$ of the SR-LASSO is not always Lipschitz continuous. 
When $\{Ax:x\in S(\bar b,\bar\lambda)\} \neq \{\bar b\}$, 
we prove (in Theorem~\ref{ca3}) that the solution mapping $S$ is locally Lipschitz continuous at $(\bar b,\bar \lambda)$ if and only if the following holds: 
\begin{condition}\label{A61}
 There exists $\bar x\in S(\bar b,\bar\lambda)$ such that $A\bar x \neq \bar b$ and $ \bar b\notin \mathrm{rge}(A_{J_{SR}})$, where
 \begin{equation}\label{J-SR}
J_{SR}=J_{SR}(\bar x):=\Big\{i\in  [n]:\Big| A_i^T\frac{\bar b-A\bar x}{\|A\bar x- \bar b\|} \Big|= \bar \lambda \Big\}.
 \end{equation}
\end{condition}
Additionally, in Theorem \ref{ca2a}, we fully characterize the local single-valued and Lipschitz continuity of $S$ around $(\bar b,\bar \lambda)$ using the following condition introduced in \cite{BBH2}: 
\begin{condition}(\cite[Assumption 2]{BBH2})\label{A2}
 For a minimizer $\bar x$ of \eqref{Pi1}, we have $A\bar x\neq \bar b$, $\ker (A_{J_{SR}})=\{0\}$, and $\bar b\notin \mathrm{rge}(A_{J_{SR}})$, where $J_{SR}$ is defined in \eqref{J-SR}.
 \end{condition}
In other words, 
we advance the sufficient result in \cite[Theorem 4.3]{BBH2} by showing that it is also a necessary condition when $\{Ax:x\in S(\bar b,\bar\lambda)\} \neq \{\bar b\}$.
These two complete characterizations of Lipschitz stability for the SR-LASSO constitute our fourth main contribution. 
Clearly, Condition \ref{A2} is equivalent to Condition \ref{A61} plus $\ker (A_{J_{SR}})=\{0\}$. 
From Proposition \ref{LCsr}, Condition \ref{A61} is equivalent to Condition \ref{conditionH1-0}.
The main contributions of this paper are summarized in Table \ref{fig1}.
\begin{table}[h]
	\caption{Sketch of the main contributions}
	\centering 
	\scriptsize
	\resizebox{1\linewidth}{!}{ 
		\begin{tabular}[c]{L{3cm} L{5cm} L{5cm}}
			\toprule 
			& Multi-valued & Single-valued \\ 
			\midrule 
			\begin{tikzpicture}[baseline=(current bounding box.center)]
				\node  (K) {General problem};
				\node  (E) [below=0.001cm of K] {\eqref{Psc}};
			\end{tikzpicture} 
			& 
			\begin{tikzpicture}[baseline=(current bounding box.center)]
				\node (A) {Condition \ref{conditionH1-0} (weak)};
				\node (A1) [rotate = 270] at ($(A.south) - (2pt, 4pt)$) {$\Rightarrow$};
				\node (B1) [anchor=west] at ($(A.south) - (-2pt, 4pt)$) {Thm. \ref{aTi1}};
				\node (B) [below=0.3cm of A] {Locally single-valued and Lipschitz of $H$};
				\node (B2) [rotate = 270] at ($(B.south) - (2pt, 4pt)$) {$\Rightarrow$};
				\node (C1) [anchor=west] at ($(B.south) - (-2pt, 4pt)$) {Thm. \ref{T01}};
				\node (C) [below=0.3cm of B] {Locally Lipschitz of $S$};
			\end{tikzpicture}
			& 
			\begin{tikzpicture}[baseline=(current bounding box.center)]
				\node (A) {Condition \ref{conditionH2-0} (strong)};
				\node (A1) [rotate = 270] at ($(A.south) - (2pt, 4pt)$) {$\Rightarrow$};
				\node (B1) [anchor=west] at ($(A.south) - (-2pt, 4pt)$) {Thm. \ref{T318-h}};
				\node (B) [below=0.3cm of A] {Locally single-valued and Lipschitz of $S$};
			\end{tikzpicture}
			\\
			\midrule 
			\begin{tikzpicture}[baseline=(current bounding box.center)]
				\node  (K) {LASSO problem};
				\node  (E) [below=0.001cm of K] {\eqref{P1}};
			\end{tikzpicture} & 
			\begin{tikzpicture}[baseline=(current bounding box.center)]
				\node (A) {No assumption};
				\node (A2) [rotate = 270] at ($(A.south) - (2pt, 4pt)$) {$\Rightarrow$};
				\node (T2) [anchor=west] at ($(A.south) - (-2pt, 4pt)$) {Thm. \ref{Ti1}};
				\node (B) [below=0.3cm of A] {Global Lipschitz of $S$};
			\end{tikzpicture}
			& 
			\begin{tikzpicture}[baseline=(current bounding box.center)]
				\node (A) {Tibshirani's sufficient condition};
				\node (A1) [rotate = 270] at ($(A.south) - (2pt, 4pt)$) {$\Leftrightarrow$};
				\node (B1) [anchor=west] at ($(A.south) - (-2pt, 4pt)$) {Thm. \ref{T318-lasso}};
				\node (B) [below=0.3cm of A] {Locally single-valued and Lipschitz of $S$};
				\node (C) [below=0.3cm of B1] {(also proved in \cite[Theorem 3.12]{nghia2024geometric})};
			\end{tikzpicture}
			\\
			\midrule 
			\begin{tikzpicture}[baseline=(current bounding box.center)]
				\node  (K) {SR-LASSO problem};
				\node  (E) [below=0.001cm of K] {\eqref{Pi1}};
			\end{tikzpicture} & 
			\begin{tikzpicture}[baseline=(current bounding box.center)]
				\node (A) {Condition \ref{A61}};
				\node (A2) [rotate = 270] at ($(A.south) - (2pt, 4pt)$) {$\Leftrightarrow$};
				\node (T2) [anchor=west] at ($(A.south) - (-2pt, 4pt)$) {Thm. \ref{ca3}};
				\node (B) [below=0.3cm of A] {Local Lipschitz of $S$};
			\end{tikzpicture}
			& 
			\begin{tikzpicture}[baseline=(current bounding box.center)]
				\node (A) {Condition \ref{A2}};
				\node (A1) [rotate = 270] at ($(A.south) - (2pt, 4pt)$) {$\Leftrightarrow$};
				\node (B1) [anchor=west] at ($(A.south) - (-2pt, 4pt)$) {Thm. \ref{ca2a}};
				\node (B) [below=0.3cm of A] {Locally single-valued and Lipschitz of $S$};
			\end{tikzpicture}
			\\
			\bottomrule 
	\end{tabular}}
	\label{fig1}
\end{table}

\textbf{Notation.} 
Denote $[n]:=\{1,2,\cdots,n\}$. For a matrix $A\in \mathbb{R}^{m\times n}$ and $I\subset [n]$, $A_I$ denotes the matrix composed of the columns of $A$ corresponding to $I$. We denote the range of $A_I$ by $\mathrm{rge}(A_I)$. 
For $x\in \mathbb{R}^n$ and $I\subset [n]$, $x_{I}$ denotes the vector in $\mathbb{R}^{|I|}$ whose entries correspond to the indices in $I$. Also, denote $I^C:=[n]\setminus I=\{i\in [n]:i\notin I\}$.
Let $\mathbb{R}_+=[0,\infty)$, $\mathbb{R}_{++}=(0,\infty)$ and $\mathbb{R}_-=(-\infty,0]$. 
For $t\in \mathbb{R}$, the sign of $t$ is given by $\mathrm{sgn}(t)$. For a vector $a\in \mathbb R^n$, denote $\mathrm{sgn}(a):=(\mathrm{sgn}(a_i))_{i=1}^n$. Its support is defined by $\mathrm{supp}(a):=\{i\in [n]:a_i\neq 0\}$.  
Let $E$ be a Euclidean space,
the Euclidean $\ell_p$-norm of a vector $x\in E$ is denoted by $\| x\|_p$. Particularly, we denote $\|x\|=\|x\|_2$.
For $\bar x\in E$ and $\delta>0$, we use $B(\bar x, \delta)$ and $B[\bar x,\delta]$ to denote the open and closed balls with center at $\bar x$ and radius $\delta$, respectively. The closed unit ball $B[0,1]$ is denoted by $\mathbb B$. 

\section{Preliminaries}
\label{section2}
\setcounter{section}{2}  \setcounter{equation}{0}

In this section, we provide some  preliminaries in variational analysis. We adopt  the notation and conventions as presented in the books by Rockafellar and Wets \cite{rockwets98}. Readers may also refer to the books by Mordukhovich \cite{M2006} and Dontchev and Rockafellar \cite{DR09}.
For a Euclidean space $E$, a closed set $A\subset E$ and a point $\bar x\in A$, the regular (Clarke) tangent cone to $A$ at $\bar x$ is defined by
\begin{equation*}
\widehat T(A,\bar x):=\left\{v\in E:\forall t^k\downarrow 0~\text{and}~x^k\stackrel{A}{\longrightarrow} \bar x~\exists z^k\stackrel{A}{\longrightarrow} \bar x~\text{with}~(z^k-x^k)/t^k\rightarrow v\right\}.
\end{equation*}
The regular (Fr\'echet) normal cone to $A$ at $\bar x$ is defined by 
\begin{equation*}
\widehat N(A,\bar x):=\{u\in E:\langle u,x-\bar x\rangle\leq o(\|x-\bar x\|)~\text{for}~ x\in A\}.
\end{equation*}
The Mordukhovich limiting normal cone to $A$ at $\bar x$ is defined by
\begin{equation*}
 N(A,\bar x):=\{u\in E: \exists u^k\rightarrow u~\text{and}~x^k\stackrel{A}{\longrightarrow} \bar x ~\text{such that}~u^k\in \widehat N(A,x^k)\}.
\end{equation*}
For a set $C\subset E$, denote the polar of $C$ by $C^-:=\{u\in E: \langle u,v\rangle\leq 0~\forall  v\in C\}. $  It is known that $\widehat T(A,\bar x)=N(A,\bar x)^-$~(see \cite[Theorem 6.28(b)]{rockwets98}). 
A set $C$ is said to be locally closed at a point $\bar x$ (not necessarily in $C$)  if
$C\cap V$ is closed for some closed neighborhood $V$ of $\bar x$. 
A set $C$ is said to be a polyhedral convex set, if there exist $c_1,\dots,c_k\in E$ and $\alpha_1,\dots,\alpha_k\in \mathbb R$ such that
$C=\{x\in E:\langle c_i,x\rangle\leq \alpha_i~\text{for}~i\in[k]\}.$  
For two sets $C$ and $D$, the excess of
$C$ beyond $D$ is defined by
$e(C,D):=\sup_{x\in C}d(x,D).$ Note that $e(C,D)\neq e(D,C).$
The Hausdorff distance between $C$ and $D$ is the quantity
$$d_H(C,D):=\max\{e(C,D),e(D,C)\}.$$
In equivalent geometric form, we have
$e(C,D)=\inf \{\tau\ge0:C\subset D+\tau \mathbb{B}\}$
and
$$d_H(C,D)=\inf \{\tau\ge0:C\subset D+\tau \mathbb{B},D\subset C+\tau \mathbb{B}\}.$$

For a lower semicontinuous convex function $f: E \rightarrow \mathbb{R}\cup \{\infty\}$, the subdifferential of $f$ at $\bar x\in \mathrm{dom}(f):=\{x\in E :f(x)<\infty\}$ is given by
\begin{equation*}
\partial f(\bar x):=\{v\in E :\langle v,x-\bar x\rangle\leq f(x)-f(\bar x)~\forall  x\in E \}.
\end{equation*}
Specifically, the subdifferentials of the $l_2$-norm and $l_1$-norm are given by
\begin{equation}\label{subdiff-L2-L1-norms}
\partial \|\cdot\|(x)=\left\{
  \aligned
 \frac{x}{\|x\|},~x\neq 0,\\
 \mathbb{B},~x=0,
  \endaligned
  \right.
  ~~\mathrm{and}~~
  \partial \|\cdot\|_1(x)=\prod\limits_{i=1}^n\left\{
  \aligned
 \mathrm{sgn}(x_i),~x_i\neq 0,\\
 [-1,1],~x_i=0.
  \endaligned
  \right.
\end{equation}
Let $F:\mathbb{R}^p\rightrightarrows \mathbb{R}^q$ be a multifunction. The domain and graph of $F$ are defined by
$$\mathrm{dom} (F):=\{x\in \mathbb{R}^p:F(x)\neq\emptyset\}
~~\mathrm{and}~~
\mathrm{gph}(F):=\{(x,y)\in \mathbb{R}^p\times \mathbb{R}^q:y\in F(x)\},$$
 respectively. Recall that $F$ is locally bounded at $\bar x\in \mathbb{R}^p$ if for some neighborhood $U$ of $\bar x$ the set $F(U)$ is bounded. It is called locally bounded on $\mathbb{R}^p$ if $F$ is locally bounded at every $\bar x\in \mathbb{R}^p$. For $\bar x\in \mathrm{dom} (F)$ and a set $D\subset \mathrm{dom}(F)$, the outer and inner limits of $F$ at $\bar x$ relative to $D$ is defined by
 $$\limsup\limits_{x\xrightarrow D \bar x}F(x):=\left\{y\in \mathbb{R}^q:\liminf \limits_{x\xrightarrow D \bar x}d(y,F(x))=0\right\},$$
 and
 $$\liminf\limits_{x\xrightarrow D \bar x}F(x):=\left\{y\in \mathbb{R}^q:\lim \limits_{x\xrightarrow D \bar x}d(y,F(x))=0\right\}.$$  
The multifunction $F$ is said to be inner semicontinuous at $\bar x\in D$ relative to $D$ if
 $F(\bar x)\subset \liminf\limits_{x\xrightarrow D \bar x}F(x)$, 
 but outer semicontinuous at $\bar x\in D$ relative to $D$ if
 $\limsup\limits_{x\xrightarrow D \bar x}F(x) \subset F(\bar x)$.  
 Additionally, $F$ is said to be inner (outer) semicontinuous relative to $D$ if $F$ is inner (outer) semicontinuous at every $\bar x\in D$ relative to $D$.
We also adopt the following coderivative $D^*F(\bar x,\bar y):\mathbb{R}^q\rightrightarrows \mathbb{R}^p$ for $F$ at $(\bar x,\bar y)\in \mathrm{gph}(F)$:
 $$D^*F(\bar x,\bar y)(v):=\left\{u\in \mathbb{R}^p:(u,-v)\in  N(\mathrm{gph}(F),(\bar x,\bar y))\right\}~~\forall v\in \mathbb{R}^q.$$ 
The multifunction $F$
is said to be polyhedral if $\mathrm{gph}(F)$ can be expressed as the union of finitely many polyhedral convex sets.
 
For a multifunction $F:\mathbb{R}^p\rightrightarrows \mathbb{R}^q$ and a set $D\subset \mathrm{dom}(F)$, $F$ is said to be outer Lipschitz continuous at $\bar x\in D$ relative to $D$ if there exist $\delta, \kappa\in(0,\infty)$ such that
$$e(F(x),F(\bar x))\leq \kappa\|x-\bar x\|~~\forall x\in B(\bar x,\delta)\cap D,$$
or in equivalent geometric form, 
$$F(x)\subset F(\bar x)+ \kappa\|x-\bar x\|\mathbb{B}~~\forall x\in B(\bar x,\delta)\cap D.$$
If $F$ is outer Lipschitz continuous at every $x\in D$ with the same $\kappa$, then $F$ is said to be outer Lipschitz continuous relative to $D$. The terminology ``outer Lipschitz continuity'' was introduced by Robinson \cite{Robinson81} under the name  ``upper Lipschitz continuity". 
For $\bar y\in F(\bar x)$, $F$ is said to have the Aubin property at $\bar x$ for $\bar y$ if $\mathrm{gph}(F)$ is locally closed at $(\bar x,\bar y)$ and there is a constant $\kappa\geq 0$ together with neighborhoods $U$ of $\bar x$ and $V$ of $\bar y$ such that
 $$F(x)\cap V \subset F(x')+\kappa \|x-x'\|\mathbb B~~\forall x, x'\in U.$$
The infimum of $\kappa$ over all such combinations of $\kappa, U, V$ is called the Lipschitz modulus of $F$ at $\bar x$ for $\bar y$ and denoted by $\mathrm{lip}\left( F; \bar x | \bar y \right)$. 
We can also define the partial Aubin property for the mapping $f:\mathbb{R}^d\times \mathbb{R}^p\rightrightarrows \mathbb{R}^q$ and denote the partial Lipschitz modulus of $f$ with respect to $p$ uniformly in $x$ at $(\bar p, \bar x)$ by $\widehat{\mathrm{lip}}_p \left(f; (\bar p, \bar x) \right)$ if $f$ is single-valued at $(\bar p, \bar x)$. For details, see \cite[Section 3E]{DR09}. 

At the end of this section, we recall the definition of Lipschitz continuity for set-valued mappings, which plays a central role in our paper. 
\begin{definition}
    A mapping $F$ is said to be (Hausdorff) Lipschitz continuous relative to $D$ if there exist $\kappa\in(0,\infty)$, a Lipschitz constant, such that
$$d_H (F(x),F(x'))\leq \kappa\|x-x'\|~~~\forall x, x'\in D,$$
or equivalently 
$F(x)\subset F(x') + \kappa\|x -x'\|\mathbb{B}$ for all $ x, x'\in D.$
$F$ is said to be locally (Hausdorff) Lipschitzian around $\bar x$ if there exist $\delta, L\in(0,\infty)$ such that
$$d_H(F(x),F(x'))\leq L\|x - x'\|~~~\forall x, x'\in B(\bar x,\delta).$$
\end{definition}

\section{The general problem}
\label{section-extended}
\setcounter{equation}{0}

This section utilizes our decomposition method from Section \ref{introduction-1} to analyze the Lipschitz stability of solution mapping:
\begin{equation*}\label{SL1}
  S:(b,\lambda)\in \mathbb{R}^m\times \mathbb{R}_{+}\mapsto \mathrm{argmin}_{x\in\mathbb{R}^n} \left\{ h(Ax-b)+\lambda \|x\|_1 \right\},
\end{equation*}
for the general LASSO-type problem \eqref{Psc}. 

\subsection{Basic properties of solution mapping}

For $(b,\lambda)\in \mathbb{R}^m\times \mathbb{R}_{++}$, since the function $x\mapsto h(Ax-b)+\lambda \|x\|_1$ is continuous and level-bounded, by \cite[Theorem 1.9]{rockwets98}, the solution set $S(b,\lambda)$ is nonempty and compact.
Furthermore, by the first-order optimality condition and classical subdifferential sum and chain rules, we have $x\in S(b,\lambda)$ if and only if $0\in A^T\partial h(Ax-b)+\lambda\partial \|\cdot\|_1(x)$. Hence, 
\begin{equation}\label{e101}
  S(b,\lambda)=\left\{x\in \mathbb{R}^n:0\in {A^T \partial h(Ax-b)}+{\lambda}\partial \|\cdot\|_1(x) \right\}.
\end{equation}

We present two results that are ubiquitous in our study. 

\begin{lemma}\label{L625a}
Let $h:\mathbb{R}^m\rightarrow \mathbb{R}_+$ be a continuous convex function. 
Then the solution mapping $S$ has a closed graph and is locally bounded at each $(b,\lambda)\in \mathbb{R}^m\times \mathbb R_{++}$.
\end{lemma}
\begin{proof}
    First, for any $(b^k,\lambda^k,x^k)\in \mathrm{gph}(S)$ with $(b^k,\lambda^k,x^k)\rightarrow (\hat b,\hat \lambda, \hat x)$, we have
    $$h(Ax^k-b^k)+\lambda^k \|x^k\|_1\leq h(Ax-b^k)+\lambda^k \|x\|_1~~\forall x\in \mathbb R^n.$$
    Letting $k\rightarrow \infty$, by the continuity of $h$, we get
    $$h(A\hat x-\hat b)+\hat \lambda \|\hat x\|_1\leq h(Ax-\hat b)+\hat\lambda \|x\|_1~~\forall x\in \mathbb R^n.$$
    This impies that $(\hat b,\hat \lambda, \hat x)\in \mathrm{gph}(S)$. Hence $\mathrm{gph}(S)$ is closed.
    
Second, since $h$ is nonnegative, we have 
$$\lambda\|x\|_1\leq h(Ax-b)+\lambda \|x\|_1\leq h(0-b)+\lambda \|0\|_1=h(-b)~~\forall x\in S(b,\lambda).$$
This implies that $S$ is locally bounded at each $( b, \lambda)$ with $\lambda>0$. 
\end{proof}
 
\begin{lemma}\label{L613a}
    Let $(b,\lambda)\in \mathbb{R}^m\times \mathbb{R}_{++}$. Suppose that $\bar x\in S(b,\lambda)$ and $h$ is continuously differentiable at $A\bar x-b$. Denote
    $$J_1:=\{i\in [n]:A_i^T\nabla h(A\bar x-b)=-\lambda\}
    ~\text{and}~
    J_2:=\{j\in [n]:A_j^T\nabla h(A\bar x-b)=\lambda\}.$$
    Define $a\in \mathbb R^n$ by
\begin{equation}\label{ea}
a_i=1~\forall i\in J_1,~a_j= -1~\forall j\in J_2,~a_k=0~\forall k\in (J_1\cup J_2)^C.
\end{equation}
Then, (i)~$\mathrm{supp}(\bar x)\subset J_1\cup J_2$; 
(ii)~For $t>0$, we have $\mathrm{sgn}(\bar x+t a)=\mathrm{sgn}(a)=a$, $\mathrm{supp}(\bar x+t a)=J_1\cup J_2$, and $\bar x+t a\in S(b+tAa,\lambda)$.
\end{lemma}
\begin{proof}
Since $\bar x$  is a solution of \eqref{Psc}, by \eqref{e101} and \eqref{subdiff-L2-L1-norms}, we have
\begin{equation}\label{key-1}
0\in A^T\nabla h(A\bar x-b)+\lambda\partial\|\cdot\|_1(\bar x)
~\mathrm{and}~
\left| A_i^T\nabla h(A\bar x-b) \right|\leq \lambda~~\forall i\in[n] .
\end{equation}
Hence $\bar x_i\ge 0$ for $i\in J_1$, $\bar x_j\leq 0$ for $j\in J_2$, $\bar x_k=0$ for $k\in (J_1\cup J_2)^C$. Then $\mathrm{supp}(\bar x)\subset J_1\cup J_2$. 
Additionally, for $t>0$, $\mathrm{sgn}(\bar x+t a)=\mathrm{sgn}(a)=a$. Thus $\mathrm{supp}(\bar x+t a)=J_1\cup J_2$. 
This implies that $\left[ \partial\|\cdot\|_1(\bar x+t a) \right]_k=[-1,1]$ for all $k\in (J_1\cup J_2)^C$. 
Hence, by \eqref{key-1} and the definition of $J_1, J_2$, we have 
$$A^T\nabla h(A(\bar x+t a)-(b+tAa))=A^T\nabla h(A\bar x-b)\in -\lambda \partial\|\cdot\|_1(\bar x+t a).$$
Therefore, $\bar x+t a\in S(b+tAa,\lambda)$. 
\end{proof}

\subsection{Single-valued and Lipschitz property of $H$}

This section is to achieve the task in Stage I of our decomposition method, which is to establish the single-valued and Lipschitz property of the mapping $H$ defined in \eqref{1521a-0} under proper conditions. 
We first prove that $H(\bar b, \bar \lambda)$  is a singleton under Condition \ref{conditionH1-0}.


\begin{proposition}\label{L514a}
Let $\bar \lambda>0$. 
Suppose that Condition \ref{conditionH1-0} holds for $(\bar b,\bar \lambda)$, that is, for some $\bar x\in S(\bar b,\bar \lambda)$, $h$ is twice continuously differentiable around $A\bar x-\bar b$ and
\begin{equation}\label{520g}
 \ker (\nabla^2h(A\bar x-\bar b))\cap \mathrm{rge}(A_{J_h})=\{0\},
\end{equation}
where $J_h=J_{h}(\bar x):=\left\{i\in [n]:|A_i^T\nabla h(A\bar x-\bar b)|=\bar\lambda\right\}$.
Then 
\begin{equation*}
   H(\bar b,\bar \lambda)= \{Ax:x\in S(\bar b,\bar \lambda)\}=\{A\bar x\}.
\end{equation*}
\end{proposition}

\begin{proof}
First, by Lemma \ref{L613a}(i) for $\bar x\in S(\bar b,\bar \lambda)$, we get $\mathrm{supp}(\bar{x})\subset J_h(\bar x)=:J_h$. Then $A\bar x \in \mathrm{rge}(A_{J_h})$. 
We claim that 
\begin{equation}\label{support-of-solutions}
    \mathrm{supp}(\hat{x})\subset J_h
    ~\mathrm{and}~\mathrm{then}~A\hat{x}\in \mathrm{rge}(A_{J_h})
    ~~\forall \hat{x}\in S(\bar b,\bar \lambda).
\end{equation}
Indeed, for any $\hat{x}\in S(\bar b,\bar \lambda)$, 
since $S(\bar b,\bar \lambda)$ is convex,  $x^k:=\frac{1}{k}\hat x+\frac{k-1}{k}\bar x\in S(\bar b,\bar \lambda)$ for all $k\geq 1$. 
Since $h$ is twice continuously differentiable at $A\bar x-\bar b$, it is also continuously differentiable at $Ax^k-\bar b$ for all sufficiently large $k$. Hence, by the definition of $J_h$, since $x^k\rightarrow \bar x$ as $k\rightarrow\infty$, for sufficiently large $k$, 
$$|A_i^T\nabla h(A x^k-\bar b)|<\bar\lambda~~\forall i\in J_h^C.$$
Thus, for such $k$, by Lemma \ref{L613a}(i), $ x^k_{J_h^C}=0$. Recall that $ \bar x_{J_h^C}=0$ and $\hat x=k x^k -(k-1) \bar{x}$, we have
$\hat x_{J_h^C}=0$. Then \eqref{support-of-solutions} holds. Thus $A\hat{x}- A\bar{x}\in \mathrm{rge}(A_{J_h})$ for all $\hat{x}\in S(\bar b,\bar \lambda)$. 

Next we will prove that 
\begin{equation}\label{520f}
  \left[\nabla^2h(A\bar x-\bar b)\right] \left[ A(\hat x-\bar x) \right]=0
  ~~\forall \hat{x}\in S(\bar b,\bar \lambda).
\end{equation}
Then, by \eqref{520g}, we get $A\hat{x}- A\bar{x}=0$, which implies that $H(\bar{b}, \bar{\lambda})$ is a singleton. 

To prove \eqref{520f}, for any $\hat{x}\in S(\bar b,\bar \lambda)$, 
since $S(\bar b,\bar \lambda)$ is convex, we have $\bar x+t(\hat x-\bar x)\in S(\bar b,\bar \lambda)$ for all $t\in [0,1]$. Thus, 
\begin{equation}\label{520d}
  h(A(\bar x+t(\hat x-\bar x))-\bar b)+\bar\lambda \|\bar x+t(\hat x-\bar x)\|_1=h(A\bar x-\bar b)+\bar\lambda \|\bar x\|_1 ~~\forall t\in [0,1].
\end{equation}
On the other hand, by the convexity of $h$ and $\|\cdot\|_1$, for all $t\in [0,1]$, we have 
 \begin{equation}\label{inequality-convex}
 \aligned
   &h(A(\bar x+t (\hat x-\bar x))-\bar b)+\bar\lambda \|\bar x+t (\hat x-\bar x)\|_1\\
   \leq & t h(A\hat x-\bar b)+(1-t )h(A\bar x-\bar b)+t \bar\lambda \|\hat x\|_1+(1-t )\bar\lambda \|\bar x\|_1\\
   =& h(A\bar x-\bar b)+\bar\lambda \|\bar x\|_1  ~~~~~(\mathrm{since}~\bar x, \hat{x} \in S(\bar b,\bar \lambda) ).
   \endaligned
 \end{equation}
Hence, by \eqref{520d}, the above inequality must be an equality. This implies that
\begin{equation}\label{520e}
h(A(\bar x+t(\hat x-\bar x))-\bar b)=th(A\hat x-\bar b)+(1-t)h(A\bar x-\bar b) ~~\forall t\in [0,1].
\end{equation}
Since $h$ is twice continuously differentiable at $A\bar x-\bar b$, by \eqref{520e}, we have 
\begin{equation}\label{520f-0}
  \left[ A(\hat x-\bar x) \right]^T \left[\nabla^2h(A\bar x-\bar b)\right] \left[ A(\hat x-\bar x) \right]=0.
\end{equation}
Because $\nabla^2h(A\bar x-\bar b)$ is symmetric and positive semidefinite due to the convexity of $h$, there is exactly one positive semidefinite and symmetric matrix $B$ such that $\nabla^2h(A\bar x-\bar b)=BB$. Substituting this into \eqref{520f-0}, we get $B\left[ A(\hat x-\bar x) \right]=0$. Then $\left[\nabla^2h(A\bar x-\bar b)\right] \left[ A(\hat x-\bar x) \right]=BB\left[ A(\hat x-\bar x) \right]=0$. 
The proposition is proved.
\end{proof}

\begin{remark}
(i) Proposition \ref{L514a} shows that the condition required for $\bar x$ in Condition \ref{conditionH1-0} holds for all $x\in S(\bar b, \bar \lambda)$ if it holds for some $\bar x \in S(\bar b, \bar \lambda)$. 
This means that Condition \ref{conditionH1-0} can be verified for any specific choice of $x\in S(\bar b, \bar \lambda)$. \\
(ii) If $h$ is strict convex, we can directly prove $A\hat{x}=A\bar{x}$ using \eqref{inequality-convex}. This is actually the proof argument used in \cite[Lemma 4.1]{ZY}. \\
(iii) When applied to the SR-LASSO, Proposition \ref{L514a} combined with Proposition \ref{LCsr} reveals that $H(b,\lambda)$ is a singleton for the SR-LASSO under Condition \ref{A61}. 
By checking the case of $\lambda=1$ in Example \ref{fla}, Condition \ref{A61} is tight, namely, $H(b,\lambda)$ may not be a singleton when Condition \ref{A61} does not hold. 
\end{remark}

Next, we will prove the representation \eqref{H-represent}, which is an important tool in establishing the Lipschitz continuity of $H$.
Let $G$ be a multifunction defined by
\begin{equation}\label{eG1}
  G(y):=\left\{z\in \mathbb R^m: \exists x\in \mathbb R^n~\text{with}~A^Tz\in \partial \|\cdot\|_1(x)~\text{and}~y=Ax \right\}
  ~~\forall y\in \mathbb R^m. 
\end{equation}

\begin{lemma}\label{L619b}
    The representation \eqref{H-represent} holds. Namely, 
\begin{equation}\label{eH1}
    H(b,\lambda)=\left\{y\in \mathbb R^m:0\in \frac{1}{\lambda} \partial h(y-b) +G(y)\right\}
    ~~\forall (b,\lambda)\in  \mathbb{R}^m\times \mathbb{R}_{++} .
\end{equation}
\end{lemma}

\begin{proof}
    By \eqref{e101}, for $\lambda>0$, the equivalent relationships below hold:
\begin{equation*}
\aligned
    y\in H(b,\lambda)\Longleftrightarrow &~\exists x\in S(b,\lambda)~\text{with}~y=Ax\\
    \Longleftrightarrow &~\exists x\in \mathbb R^n~\text{with}~y= Ax~ \text{and} ~0\in \frac{1}{\lambda} A^T\partial h(Ax-b) +\partial \|\cdot\|_1(x)\\
    \Longleftrightarrow &~\exists z\in \frac{\partial h(y-b)}{\lambda}~\text{and}~x\in \mathbb R^n~\text{with}~y= Ax~ \text{and}~ A^T(-z)\in \partial \|\cdot\|_1(x)\\
    \Longleftrightarrow &~\exists z\in \frac{\partial h(y-b)}{\lambda}~\text{with}
    ~ -z\in G(y)\\
    \Longleftrightarrow &~ 
    0\in \frac{1}{\lambda} \partial h(y-b) +G(y) .
    \endaligned
\end{equation*}
This implies the representation \eqref{eH1}.
\end{proof}

\begin{lemma}\label{Lmon}
The inverse mapping $G^{-1}$ of $G$ satisfies 
\begin{equation}\label{inverse-G}
    G^{-1}(z)= A(\partial \|\cdot\|_1)^{-1}(A^Tz)
    ~~\forall z\in \mathbb R^m.
\end{equation}
And both $G$ and $G^{-1}$ are maximal monotone.
\end{lemma}

\begin{proof}
By the definition of $G$ in \eqref{eG1}, we have 
\begin{equation*}
\aligned
    z\in G(y)\iff & \exists x\in \mathbb R^n~\text{such that}~A^Tz\in \partial \|\cdot\|_1(x)~\text{and}~y=Ax \\
    \iff &\exists x\in \mathbb R^n~\text{such that}~y=Ax~ \text{and} ~x\in (\partial \|\cdot\|_1)^{-1}(A^Tz)\\
    \iff & y\in A(\partial \|\cdot\|_1)^{-1}(A^Tz).
    \endaligned
\end{equation*}
This proves \eqref{inverse-G}. 
By \cite[Theorem 12.17 and Exerise 12.8(a)]{rockwets98}, both $\partial \|\cdot\|_1$ and $(\partial \|\cdot\|_1)^{-1}$ are maximal monotone. 
Since $0\in\mathrm{rint}(\mathrm{dom}((\partial \|\cdot\|_1)^{-1}))$, by \eqref{inverse-G} and \cite[Theorem 12.43 and Exerise 12.8(a)]{rockwets98}, both $G^{-1}$ and $G$ are maximal monotone. 
\end{proof}

Now, we can use the representation \eqref{H-represent} and the fact that $G$ is a maximal monotone mapping to establish the single-valued and Lipschitz continuity of $H$ under Condition \ref{conditionH1-0}. 
Before doing this, we recall a notable coderivative criterion for generalized equations from \cite[Exercise 4C.5]{DR09}. 
A set-valued mapping $\Phi:\mathbb R^r\rightrightarrows \mathbb R^s$ is said to be positively
homogeneous if $\mathrm{gph}(\Phi)$ is a cone, and its outer norm is defined by
$|\Phi|^+:=\sup\limits_{\|w\|\leq 1}\sup\limits_{z\in \Phi(w)}\|z\|.$ 

\begin{lemma}(\cite[Exercise 4C.5]{DR09})\label{4C.5}
For the solution mapping $\Psi(p):=\{x :0\in f(p,x)+F(x)\}$ and a pair $(\bar p,\bar x)$ with $\bar x\in \Psi(\bar p)$, suppose that $f$ is differentiable in a neighborhood of $(\bar p,\bar x)$, $\mathrm{gph}(F)$ is locally closed at $(\bar x,-f(\bar p,\bar x))$, and 
$$|\left(D_x^*f(\bar p,\bar x)+D^*F(\bar x,-f(\bar p,\bar x))\right)^{-1}|^+\leq \mu <\infty.$$
Then $\Psi$ has the Aubin property at $\bar p$ for $\bar x$ with 
\begin{equation}\label{lipcon}
    \mathrm{lip}\left( \Psi ; \bar p | \bar x \right) 
    \leq \mu\, \widehat{\mathrm{lip}}_p
    \left(f; (\bar p, \bar x) \right).
\end{equation}
\end{lemma}

\begin{theorem}\label{aTi1}
Let $h:\mathbb{R}^m\rightarrow \mathbb{R}_+$ be a continuous convex function. Given $(\bar b,\bar\lambda)\in \mathbb R^m\times \mathbb R_{++}$, suppose that Condition \ref{conditionH1-0} holds. Then the mapping $H$ defined in \eqref{1521a-0} is single-valued and Lipschitz continuous in a neighborhood of $(\bar b,\bar \lambda)$. 
\end{theorem}

\begin{proof}
We divide the proof into two steps.

{\bf \textit{Step 1. Local single-valuedness of $H(b,\lambda)$ around $(\bar b,\bar \lambda)$.} }

First, Condition \ref{conditionH1-0} implies that there exists $\bar x\in S(\bar b,\bar \lambda)$ such that $h$ is twice continuously differentiable around $A\bar x-\bar b$ and $\ker(\nabla^2h(A\bar x- \bar b))\cap \mathrm{rge}(A_{J_h})=\{0\}$, where $J_h=\left\{i\in [n]:|A_i^T\nabla h(A\bar x-\bar b)|=\bar\lambda\right\}$. Thus we have $\nabla^2h(A\bar x- \bar b)(z)\neq 0$ for all $z\in \mathrm{rge}(A_{J_h})$ satisfying $\| z \|=1$. Hence, there exists $\gamma>0$ such that 
\begin{equation}\label{+619a}
    \|\nabla^2h(A\bar x- \bar b)(z)\|\ge\gamma>0 ~~\forall z\in \mathrm{rge}(A_{J_h})~\mathrm{satisfying}~\| z \|=1  .
\end{equation}
By the definition of $J_h$ and \eqref{e101}, 
   $ |A_i^T\nabla h(A\bar x-\bar b)|<\bar \lambda$ for all $ i\in J_h^C.$
Hence, 
there exists $\delta_1\in(0,\bar \lambda)$ such that for all $(y,b,\lambda)\in B(A\bar x,\delta_1)\times B(\bar b,\delta_1)\times B(\bar\lambda,\delta_1)$, 
\begin{equation}\label{e619b}
   |A_i^T\nabla h(y- b)|< \lambda~~\forall i\in J_h^C,
\end{equation}
and $\|\nabla^2h(y-  b)(z)\|>0$ for all $z\in \mathrm{rge}(A_{J_h})$ satisfying $\| z \|=1$. 
This implies that 
\begin{equation}\label{+619b}
    \ker(\nabla^2h(y-  b))\cap \mathrm{rge}(A_{J_h})=\{0\}~~\forall(y,b,\lambda)\in B(A\bar x,\delta_1)\times B(\bar b,\delta_1)\times B(\bar\lambda,\delta_1).
\end{equation}

Second, by Lemma \ref{L625a}, $S$ is locally bounded at $( \bar b,\bar \lambda)$. Consequently, $H$ is locally bounded at $(\bar b,\bar \lambda)$. 
Since $\mathrm{gph}(H)$ is closed by Lemma \ref{L625a}, $H$ is outer semicontinuous at $(\bar b,\bar\lambda)\in \mathrm{dom}(H)$. Given that $H(\bar b,\bar \lambda)=\{ A\bar x \}$ due to Proposition \ref{L514a}, by \cite[Theorem 5.19]{rockwets98}, there exists $\delta_2\in(0,\delta_1)$ such that $H(b,\lambda)\subset B(A\bar x,\delta_1)$~for all $(b,\lambda)\in B(\bar b,\delta_2)\times B(\bar\lambda,\delta_2)$.
Hence, for such $(b,\lambda)$, we have $Ax\in B(A\bar x,\delta_1)$ for all $x\in S(b,\lambda)$. 
Thus, by \eqref{e619b} and \eqref{+619b}, we have 
 \begin{equation}\label{510e}
 \max\limits_{i\in J_h^C}|A_i^T\nabla h(Ax- b)|< \lambda~\text{and}~\ker(\nabla^2h(Ax-  b))\cap \mathrm{rge}(A_{J_h})=\{0\} ~~\forall x\in S(b,\lambda).
 \end{equation}
This implies that for all $(b,\lambda)\in  B(\bar b,\delta_2)\times B(\bar \lambda,\delta_2)$ and $x\in S(b,\lambda)$, 
 $$J_{h(b,\lambda)}:=\{i\in [n]:|A_i^T\nabla h(Ax- b)|= \lambda\}\subset J_h .$$
Therefore, $\ker(\nabla^2h(A x-  b))\cap  \mathrm{rge}(A_{J_{h(b,\lambda)}})\subset \ker(\nabla^2h(A x-  b))\cap \mathrm{rge}(A_{J_h}) =\{0\}$. 
Using Proposition \ref{L514a}, $H( b,\lambda)$ is single-valued for all $(b,\lambda)\in  B(\bar b,\delta_2)\times B(\bar \lambda,\delta_2)$.

{\bf\textit{Step 2. 
Local Lipschitz continuity of $H(b,\lambda)$ around $(\bar b,\bar \lambda)$.} }

Define a mapping $f: \mathbb R^m\times \mathbb R_{++}\times \mathbb R^m \rightrightarrows \mathbb R^m$ by 
\begin{equation*}
f(b,\lambda, y):=\frac{1}{\lambda } \partial h(y-b) ~~ \forall (b,\lambda, y)\in \mathbb R^m\times \mathbb R_{++}\times \mathbb R^m.
\end{equation*}
By Lemma \ref{L619b}, we have 
\begin{equation}\label{e619cc}
     H(b,\lambda)=\left\{y\in \mathbb R^m:0\in f(b,\lambda, y) +G(y)\right\}~~\forall (b,\lambda)\in \mathbb R^m\times \mathbb R_{++}.
\end{equation}
Since $h$ is twice continuously differentiable at $A\bar x- \bar b$, $f$ is differentiable in a neighborhood of $(\bar b,\bar \lambda, A\bar x)$. 
Denote $\bar y:=A\bar x$, 
using $\delta_2$ from Step 1, for all $(b,\lambda, y)\in B(\bar b,\delta_2)\times B(\bar\lambda,\delta_2)\times B(\bar y,\delta_2)$, we have  
\begin{equation*}
\begin{aligned}
f(b,\lambda, y)=\frac{1}{\lambda } \nabla h(y-b)  ~\mathrm{and}~
D_y^*f(\bar b, \bar \lambda, \bar y)=\frac{1}{\bar \lambda } \nabla^2 h(A\bar x- \bar b) .
\end{aligned}
\end{equation*}

For any $v\in \mathbb{R}^m $ with $\| v\|\leq1$ and
   $u\in (D_y^*f(\bar b,\bar\lambda, \bar y)+D^*G(\bar y,-f(\bar b,\bar\lambda, \bar y)))^{-1}(v),$
we have
$v\in D_y^*f(\bar b,\bar\lambda, \bar y)(u)+D^*G(\bar y,-f(\bar b,\bar\lambda, \bar y))(u)$. 
Thus, 
\begin{equation}\label{relation-u-v-0}
v-D_y^*f(\bar b,\bar \lambda, \bar y)(u)\in D^*G(\bar y,-f(\bar b,\bar\lambda, \bar y))(u).
\end{equation}
By \cite[Theorem 5.6]{M2018} and Lemma \ref{Lmon}, it follows that
\begin{equation}\label{dd}
    \left\langle v-D_y^*f(\bar b,\bar \lambda, \bar y)(u), u \right\rangle \geq 0.
\end{equation}
Note that $D_y^*f(\bar b,\bar \lambda, \bar y)$ is symmetric and positive semidefinite, and by Condition \ref{conditionH1-0}, 
$$\ker(D_y^*f(\bar b,\bar \lambda, \bar y))\cap \mathrm{rge}(A_{J_{h}})=\ker(\nabla^2 h(A\bar x-\bar b) )\cap \mathrm{rge}(A_{J_{h}})=\{0\}.$$ 
 Hence, there exists $\mu\geq 0$ such that
 \begin{equation}\label{e72}
     \|y\|^2\leq \mu 
     \left\langle D_y^*f(\bar b,\bar \lambda, \bar y) (y), y \right\rangle ~~~\forall y\in \mathrm{rge}(A_{J_{h}}).
 \end{equation}

Note that $u\in \mathrm{dom}(D^*G(\bar y,-f(\bar b,\bar\lambda, \bar y)))$ by \eqref{relation-u-v-0}. We claim that $u\in \mathrm{rge}(A_{J_{h}})$. By the definition of $D^* G$, it suffices to prove that 
\begin{equation}\label{dom-G-u-A}
    N(\mathrm{gph}(G),(\bar y,-f(\bar b,\bar\lambda, \bar y)))\subset 
    \mathbb R^m\times \mathrm{rge}(A_{J_{h}}).
\end{equation}
Recall that $\widehat T(\mathrm{gph}(G),(\bar y,-f(\bar b,\bar\lambda, \bar y)))=N(\mathrm{gph}(G),(\bar y,-f(\bar b,\bar\lambda, \bar y)))^-$ (by \cite[Theorem 6.28(b)]{rockwets98}), and $\mathbb R^m\times \mathrm{rge}(A_{J_{h}})=(\{0\}\times \ker(A_{J_h}^T))^-$. It remains to prove 
 \begin{equation}\label{e72+d}
     \{0\}\times \ker(A_{J_h}^T)\subset \widehat T(\mathrm{gph}(G),(\bar y,-f(\bar b,\bar\lambda, \bar y))) .
 \end{equation}
Indeed, let $w\in \ker(A_{J_h}^T)$. For any sequences $t^k\downarrow 0$ and $(y^k,z^k)\stackrel{\mathrm{gph}(G)}{\longrightarrow} (\bar y,-f(\bar b,\bar\lambda, \bar y))$, by the definition \eqref{eG1} of $G$, there exists $x^k\in \mathbb R^n$ such that 
 \begin{equation}\label{e72+c}
     A^Tz^k\in \partial \|\cdot\|_1(x^k)~\text{and}~y^k=Ax^k.
 \end{equation}
We claim that $A^T(z^k+t^kw)\in \partial \|\cdot\|_1(x^k)$ for all sufficiently large $k$. Then, by \eqref{e72+c}, $z^k+t^kw\in G(y^k)$. Hence $(y^k,z^k+t^kw)\stackrel{\mathrm{gph}(G)}{\longrightarrow} (\bar y,-f(\bar b,\bar\lambda, \bar y))$.  Since $[(y^k,z^k+t^kw)-(y^k,z^k)]/t^k \rightarrow (0,w)$, $(0,w)\in\widehat T(\mathrm{gph}(G),(\bar y,-f(\bar b,\bar\lambda, \bar y)))$. Thus \eqref{e72+d} holds. 
 
 To prove $A^T(z^k+t^kw)\in \partial \|\cdot\|_1(x^k)$ for all sufficiently large $k$, note that
 \begin{equation}\label{e72+a}
     A_i^T(-f(\bar b,\bar\lambda, \bar y))=-\frac{1}{\bar\lambda }A_i^T \nabla h(A\bar x-\bar b)\in (-1,1)~~\forall i\in J_h^C.
 \end{equation}
Since $z^k\rightarrow -f(\bar b,\bar\lambda, \bar y)$ and $z^k+t^kw\rightarrow -f(\bar b,\bar\lambda, \bar y)$,
there exists $k_0>0$ such that
 \begin{equation}\label{e72+bb}
     A_i^Tz^k\in (-1,1)~\text{and}~A_i^T(z^k+t^kw)\in (-1,1)~~\forall i\in J_h^C,~\forall k\geq k_0.
 \end{equation}
Hence, since $A^Tz^k\in \partial \|\cdot\|_1(x^k)$, we have $x^k_i=0$ for all $i\in J_h^C$. Thus, for each $i\in J_h^C$, $(\partial \|\cdot\|_1(x^k))_i=[-1,1]$. Then $A_i^T(z^k+t^kw)\in (\partial \|\cdot\|_1(x^k))_i$ by \eqref{e72+bb}. 
On the other hand, since $w\in \ker(A_{J_h}^T)$,  $A_j^T(z^k+t^kw)=A_j^Tz^k\in \left[\partial \|\cdot\|_1(x^k)\right]_j$ for all $j\in J_h$. This implies that $A^T(z^k+t^kw)\in \partial \|\cdot\|_1(x^k)$. Hence, we have proven $u\in \mathrm{rge}(A_{J_{h}})$. 
 
 Since $u\in \mathrm{rge}(A_{J_{h}})$, by \eqref{e72} and \eqref{dd}, we have
 \begin{equation}\label{kerDQ}
   \|u\|^2
   \leq \mu 
   \left\langle D_y^*f(\bar b,\bar \lambda, \bar y) (u), u \right\rangle
   \leq 
   \mu \langle v, u\rangle\leq \mu \|v\|\|u\|.
 \end{equation}
 This implies that $\|u\|\leq \mu \|v\|\leq \mu$. Hence, we have 
 \begin{equation}\label{lipH}
     \left|\left(D_y^*f(\bar b,\bar \lambda, \bar y)+D^*G(\bar y,-f(\bar b,\bar\lambda, \bar y))\right)^{-1}\right|\leq \mu<\infty.
 \end{equation}
 By Lemma \ref{4C.5} for \eqref{e619cc}, $H$ has the Aubin property at $(\bar b,\bar\lambda)$ for $\bar y$. That is, there exist $L_H>0$ and $\delta_3\in (0,\delta_2)$ such that
$$ H(b,\lambda)\cap B(\bar y,\delta_3)\subset H(b',\lambda')+ L_H \|(b,\lambda)-(b',\lambda')\| \mathbb{B}~$$
for all $ (b,\lambda),~(b',\lambda')\in B(\bar b,\delta_3)\times B(\bar \lambda ,\delta_3).$ 
Recalling that $H$ is out semicontinuity at $(\bar b,\bar\lambda)$, by \cite[Theorem 5.19]{rockwets98}, there exists $\delta_4\in(0,\delta_3)$ such that $ H(b,\lambda)\subset B(\bar y,\delta_3)$ for all $ (b,\lambda)\in B(\bar b,\delta_4)\times B(\bar \lambda ,\delta_4).$ 
This implies that
$$ H(b,\lambda)\subset H(b',\lambda')+ L_H \|(b,\lambda)-(b',\lambda')\| \mathbb{B}~~\forall (b,\lambda), (b',\lambda')\in B(\bar b,\delta_4)\times B(\bar \lambda ,\delta_4).$$
Hence $H$ is Lipschitz continuous on $B(\bar b,\delta_4)\times B(\bar \lambda ,\delta_4)$. 
\end{proof}

\begin{remark}
(i) Using the representation \eqref{eH1} of $H$ and a notable result from \cite{BBH}, 
we can provide another proof of the local single-valued and Lipschitz property of $H$. 
Indeed, define the mapping $Q(\cdot):=f(\bar b,\bar \lambda,\cdot)+G(\cdot)$. By \eqref{kerDQ} for $v=0$, we have 
 $$\ker D^*Q(\bar y,0)=\ker \left(D_y^*f(\bar b,\bar \lambda, \bar y)+D^*G(\bar y,-f(\bar b,\bar\lambda, \bar y))\right)=\{0\}.$$ Note that $f(\bar b,\bar \lambda,\cdot)$ and $G$ are monotone.
 Applying \cite[Proposition 4.11(b)]{BBH}, the mapping $H$ is locally single-valued and Lipschitz at $(\bar b,\bar\lambda)$. \\
 (ii) The argument in Step 1 actually shows that Condition \ref{conditionH1-0} is a local property, namely, stable to small pertubations of the data. 
\end{remark}

\subsection{Lipschitz stability of solution mapping}

This section is devoted to Stage II of our decomposition method. 
The key is to find a way to recover $S$ from $H$ that allows Lipschitz continuity to transfer from $H$ to $S$. 

Our main idea is to decompose $S$ into the union of finitely many simpler set-valued mappings.
Specifically, 
for $J_1,J_2\subset  [n]$ with $J_1\cap J_2=\emptyset$, define two simple polyhedral sets: 
\begin{equation}\label{ePE}
\aligned
  P_{J_1J_2}:=&
  \left\{ 
  x\in\mathbb{R}^n: x_i\geq 0~\forall i\in J_1,~x_i\leq 0~\forall i\in J_2,~x_i=0~\mathrm{otherwise}
  \right\}, \\
  E_{J_1 J_2}:=&
  \left\{ 
  x\in\mathbb{R}^n: x_i=1~\forall i\in J_1,~x_i=-1~\forall i\in J_2,~x_i\in[-1,1]~\mathrm{otherwise}
  \right\},
  \endaligned
\end{equation}
and a set-valued mapping $S_{J_1 J_2}:\mathbb{R}^m\times \mathbb{R}_+\rightrightarrows \mathbb{R}^n$,
\begin{equation}\label{510g}
  S_{J_1 J_2}(b,\lambda):=\left\{x\in P_{J_1 J_2}: 0\in  {A^T\partial h(Ax-b)}+ {\lambda} E_{J_1 J_2}\right\}.
\end{equation}
Clearly, $S_{J_1 J_2}(b,\lambda)$ is simpler than $S(b,\lambda)$ because $\partial \|\cdot\|_1(x)$ in $S$ is replaced by a simpler polyhedral set $E_{J_1 J_2}$ independent of $x$. 

One important observation in our study is that the solution mapping $S$ in \eqref{e101} can be decomposed by $S_{J_1 J_2}(b,\lambda)$'s. 

\begin{lemma}\label{sfj}
Let $h:\mathbb{R}^m\rightarrow \mathbb{R}_+$ be a convex function
Then, 
\begin{align}
S(b,\lambda)=&\bigcup\limits_{J_1,J_2\subset  [n],~J_1\cap J_2=\emptyset}S_{J_1 J_2}(b,\lambda)~~\forall (b,\lambda)\in \mathbb{R}^m\times \mathbb{R}_+  .
\label{510f} 
\end{align}
\end{lemma}

\begin{proof}
For $x\in\mathbb{R}^n$, define \( J_1(x) := \{i \in [n] : x_i > 0\} \) and \( J_2(x) := \{i \in [n] : x_i < 0\} \). When \( x \in P_{J_1 J_2} \), \( J_1(x) \subset J_1 \), \( J_2(x) \subset J_2 \). Thus, $E_{J_1 J_2}\subset \partial\|\cdot\|_1(x)$. 
By \eqref{e101} and \eqref{510g}, we conclude that \( S_{J_1 J_2}(b,\lambda) \subset S(b,\lambda) \) for all \( J_1, J_2 \subset [n] \) with $J_1\cap J_2=\emptyset$. 

Conversely, for any \( x \in S(b,\lambda) \), we have \( x \in P_{J_1(x) J_2(x)} \) and
\[ 0\in A^T\partial h(Ax-b)+{\lambda}\partial \|\cdot\|_1(x)={A^T\partial h(Ax-b)}+ {\lambda} E_{J_1(x) J_2(x)}. \]
This implies \( x \in S_{J_1(x) J_2(x)}(b,\lambda) \). Hence \eqref{510f} is proved.
\end{proof}

The following result shows that Lipschitz continuity can be transferred from $H$ to $S_{J_1 J_2}$. 
To this end, the key point is that, when $H(b,\lambda)$ is single-valued on a set $U\subset  \mathbb{R}^m\times \mathbb{R}_{++}$, 
the relationship \eqref{S-J1J2-F} holds, that is, 
\begin{equation}\label{624b}
     S_{J_1 J_2}(b,\lambda)=\left( F_{J_1 J_2}\circ H \right) (b,\lambda)
     ~~\forall (b,\lambda)\in \mathrm{dom}(S_{J_1 J_2})\cap U,
\end{equation}
where $F_{J_1 J_2}:\mathbb R^m\rightrightarrows \mathbb{R}^n$ is a convex polyhedral multifunction defined by 
$$F_{J_1 J_2}(y)=\{x\in P_{J_1 J_2}:Ax=y\}~~\forall y\in \mathbb R^m.$$
Note that by \cite[Theorem 2.207]{Shapiro00}, there is a constant $\kappa_{J_1 J_2}$ such that $F_{J_1 J_2}$ is Lipschitz continuous on $\mathrm{dom}(F_{J_1 J_2})$ with modulus $\kappa_{J_1 J_2}$.

\begin{lemma}\label{lemma-S-lip}
(i) The set $\mathrm{dom}(S_{J_1 J_2})$ is locally closed at each $(b, \lambda)$ with $\lambda > 0$.
\\
(ii) Suppose that $H$ is single-valued on $U\subset  \mathbb{R}^m\times \mathbb{R}_{++}$. 
Then, \eqref{624b} holds for all $J_1,J_2\subset  [n]$ with $J_1\cap J_2=\emptyset$. As a result,  
$S_{J_1 J_2}$ is Lipschitz continuous on $\mathrm{dom}(S_{J_1 J_2})\cap  U$ with modulus $\kappa_{J_1 J_2} L_H$ if $H$ is further Lipschitz continuous on $U$ with modulus $L_H$, where $\kappa_{J_1 J_2}$ is the Lipschitz modulus of $F_{J_1 J_2}$.
\end{lemma}

\begin{proof}
(i) For any $(b^k,\lambda^k)\in \mathrm{dom}(S_{J_1 J_2})$ satisfying $(b^k,\lambda^k)\rightarrow (\hat b,\hat \lambda)$ with $\hat \lambda>0$,
take $x^k\in S_{J_1 J_2}(b^k,\lambda^k)$ for all $k\in \mathbb{N}$. By the definition \eqref{510g} of $S_{J_1 J_2}$, there exists $z^k\in \partial h(Ax^k-b^k)$ such that 
\begin{equation}\label{625a}
    -A^Tz^k/\lambda^k \in  E_{J_1 J_2}.
\end{equation}
Since $\hat \lambda>0$, by Lemma \ref{L625a}, $\{x^k\}_{k=1}^\infty$ is bounded. Hence, we can take a subsequence of $\{x^k\}_{k=1}^\infty$, still denoted by $\{x^k\}_{k=1}^\infty$, such that it converges to some $\hat{x} \in \mathbb{R}^n$. 
Since $\partial h$ is locally bounded, $\{z^k\}_{k=1}^\infty$ is also bounded. We can assume that $\{z^k\}_{k=1}^\infty$ converge to some $\hat{z}\in \mathbb R^m$. 
Since $\hat \lambda>0$ and $E_{J_1 J_2}$ is closed, letting $k\rightarrow \infty$ in \eqref{625a}, we have 
$- A^T {\hat{z}} /{\hat \lambda}\in E_{J_1 J_2}$. This implies that $-A^T\hat{z}\in\hat \lambda  E_{J_1 J_2}$. Noting that $\mathrm{gph}(\partial h(\cdot))$ is closed, we have  $\hat{z}\in \partial h(A\hat x-\hat b)$. Thus $0\in A^T\partial h(A\hat x-\hat b)+\hat \lambda  E_{J_1 J_2}$. Hence, $(\hat b,\hat \lambda)\in  \mathrm{dom}(S_{J_1 J_2})$. This implies that $\mathrm{dom}(S_{J_1 J_2})$ is closed in $\mathbb{R}^m\times\mathbb{R}_{++}$.
As a result, the set $\mathrm{dom}(S_{J_1 J_2})$ is locally closed at each $(b, \lambda)$ with $\lambda > 0$.

(ii) First, note that $P_{J_1 J_2}$ is a convex polyhedral set. 
It is easy to verify that $$\mathrm{gph}(F_{J_1 J_2}):= \{(y,x)\in \mathbb{R}^m \times P_{J_1 J_2}:y=Ax\} $$ is also a convex polyhedral set. By \cite[Theorem 2.207]{Shapiro00}, there is $\kappa_{J_1 J_2}>0$ such that
\begin{equation}\label{512a}
  d_H(F_{J_1 J_2}(y),F_{J_1 J_2}(y'))\leq \kappa_{J_1 J_2} \|y-y'\|~~~\forall y, y'\in \mathrm{dom}(F_{J_1 J_2}).
\end{equation}

Second, when $H(b,\lambda)$ is single-valued on $U$,  
we prove \eqref{624b} by verifying
\begin{equation}\label{520a}
\aligned
  S_{J_1 J_2}(b,\lambda)
  =\{x\in P_{J_1 J_2}:Ax \in H(b,\lambda)\}~~\forall (b,\lambda)\in \mathrm{dom}(S_{J_1 J_2})\cap U.
  \endaligned
\end{equation}
Indeed, the ``$\subset$" of \eqref{520a} comes from \eqref{510f} and the definitions of $S_{J_1 J_2}$ and $H$.
Conversely, for any $(b,\lambda)\in \mathrm{dom}(S_{J_1 J_2})\cap U$ and $x\in P_{J_1 J_2}$ with $Ax \in H(b,\lambda)$, since $ S_{J_1 J_2}(b,\lambda)\neq\emptyset$ and $H(b,\lambda)$ is single-valued,  we have $Ax=Ax_0$ for some $x_0\in S_{J_1 J_2}(b,\lambda)$. 
Hence, we have
$$0\in  {A^T\partial h(Ax_0-b)}+ {\lambda} E_{J_1 J_2}= {A^T\partial h(Ax-b)}+{\lambda}  E_{J_1 J_2}.$$
This implies that $x\in S_{J_1 J_2}(b,\lambda)$ (by \eqref{510g}). We obtain the ``$\supset$" in \eqref{520a}. 

Third, by the assumption that
    $H$ is single-valued and Lipschitz continuous on $U$ with modulus $L_H$, we have 
\begin{equation}\label{510h}
  \|H(b,\lambda)-H(b',\lambda')\|\leq L_H \|(b,\lambda)-(b',\lambda')\|~~\forall (b,\lambda), (b',\lambda')\in U.
\end{equation}
It follows from \eqref{624b} and \eqref{512a}  that
\begin{equation*}
    \aligned
d_H(S_{J_1 J_2}(b,\lambda),S_{J_1 J_2}(b',\lambda'))
=&d_H(F_{J_1 J_2}(H(b,\lambda)),F_{J_1 J_2}(H(b',\lambda')))\\
\leq &\kappa_{J_1 J_2} \|H(b,\lambda)-H(b',\lambda')\|\\
\leq & \kappa_{J_1 J_2} L_H \|(b,\lambda)-(b',\lambda')\|
~~\forall (b,\lambda),~(b',\lambda')\in U.
    \endaligned
\end{equation*}
Hence, $S_{J_1 J_2}$ is Lipschitzian on $\mathrm{dom}(S_{J_1 J_2})\cap  U$ with modulus $\kappa_{J_1 J_2} L_H$.
\end{proof}
\begin{remark}
    It is important to note that the relationship \eqref{624b} may not hold for all  $(b,\lambda)\in U$ because $\mathrm{dom}(F_{J_1 J_2}\circ H)$ can be strictly larger than $\mathrm{dom}(S_{J_1 J_2})$. For instance, take $J_1=\{1\}$, $J_2=\{2,3\}$ for the LASSO in Example \ref{fla}, then $\mathrm{dom}(S_{J_1 J_2})=\emptyset$ but $\mathrm{dom}(F_{J_1 J_2}\circ H)\neq\emptyset$. 
    This prevents us from combining the decomposition \eqref{510f} and the relationship \eqref{624b} to directly recover $S$ from $H$.
\end{remark}

Although $S$ is expressed as the union of finitely many $S_{J_1 J_2}$, and the Lipschitz continuity of $S_{J_1 J_2}$ is established under proper conditions, we cannot directly derive the Lipschitz continuity of $S$. However, we can obtain the outer Lipschitz continuity from the decomposition instead.

\begin{lemma}\label{L516}
Let $\Phi, \Phi_1,\cdots, \Phi_l:\mathbb R^p \rightrightarrows\mathbb R^q$ be multifunctions such that $\mathrm{gph}(\Phi)=\bigcup\limits_{i=1}^l\mathrm{gph}(\Phi_i)$, and let $U\subset\mathbb R^p$. 
For all $i\in[l]$, suppose that $\Phi_i$ is Lipschitz continuous on $U\cap \mathrm{dom}(\Phi_i)$ with modulus $\kappa_i$, and $\mathrm{dom}(\Phi_i)$ is locally closed at each $z\in U$. Then $\Phi$ is outer Lipschitz continuous on $U\cap \mathrm{dom}(\Phi)$ with modulus $\max\limits_{i\in [l]}\kappa_i$.
\end{lemma}

\begin{proof}
For any $\bar z\in U\cap \mathrm{dom}(\Phi)$, denote $I(\bar z):=\{i\in[l]: \bar z\in \mathrm{dom}(\Phi_i)\}$.
Since each $\mathrm{dom}(\Phi_i)$ is locally closed at $\bar z$,  we have $d(\bar z, \mathrm{dom}(\Phi_i))>0$ if $\bar z\notin \mathrm{dom}(\Phi_i)$. Thus 
$\bar\alpha:=d(\bar z, \bigcup\limits_{i\notin I(\bar z)}\mathrm{dom}(\Phi_i) )>0$. 
Hence, for all $z\in B(\bar z, \bar\alpha)\cap \mathrm{dom}(\Phi)$, we get $I(z)\subset I(\bar z)$. This implies that $\Phi_i(\bar z)\neq\emptyset$ for each $i\in I(z)$. Then
\begin{equation*}
  \aligned
  e(\Phi(z),\Phi(\bar z))=&e(\bigcup\limits_{i\in I(z)}\Phi_i(z), \bigcup\limits_{i\in I(\bar z)}\Phi_i(\bar z))
  \leq e(\bigcup\limits_{i\in I(z)} \Phi_i(z), \bigcup\limits_{i\in I(z)}\Phi_i(\bar z))\\
 &\leq \max\limits_{i\in I(z)} 
 e(\Phi_i(z), \Phi_i(\bar z) )
  \leq \max\limits_{i\in I(z)}\kappa_i\|z-\bar z\|, 
  \endaligned
\end{equation*}
for all $z\in B(\bar z, \bar \alpha)\cap U\cap \mathrm{dom}(\Phi)$. Hence, $\Phi$ is outer Lipschitzian on $U\cap \mathrm{dom}(\Phi)$.  
\end{proof}

\begin{remark}
    Here is an example showing that we cannot directly derive the Lipschitz continuity of  $\Phi$  from the Lipschitz continuity of  $\{\Phi_i\}_{i=1}^\ell$.
   Define $\Phi_1, \Phi_2$ and $\Phi:\mathbb R\rightrightarrows \mathbb R$ by $\mathrm{gph}(\Phi_1)=(-\infty,0]\times \mathbb R$,~$\mathrm{gph}(\Phi_2)=[0,\infty)\times \mathbb R_+$ and $\mathrm{gph}(\Phi)=\mathrm{gph}(\Phi_1)\cup \mathrm{gph}(\Phi_2) $. It is easy to see that $\Phi_i$ is Lipschitz continuous on $\mathrm{dom}(\Phi_i)$ for $i=1,2$. However, $\Phi$ is not Lipschitz continuous on $\mathbb R=\mathrm{dom}(\Phi)$.
\end{remark}

\begin{corollary}\label{S-outer Lipschitzian}
Suppose that $H$ is single-valued and Lipschitz continuous on $U\subset  \mathbb{R}^m\times \mathbb{R}_{++}$ with modulus $L_H$. 
Then $S$ is outer Lipschitzian on $U$ with modulus $\max\limits_{J_1,J_2\subset  [n],~J_1\cap J_2=\emptyset} \kappa_{J_1 J_2} L_H$.
\end{corollary}
\begin{proof}
By Lemma \ref{lemma-S-lip}, for all $J_1,J_2\subset  [n]$ with $J_1\cap J_2=\emptyset$, $S_{J_1 J_2}$ is Lipschitzian on $\mathrm{dom}(S_{J_1 J_2})\cap  U$ with modulus $\kappa_{J_1 J_2} L_H$, and $\mathrm{dom}(S_{J_1 J_2})$ is locally closed at each $(b,\lambda)\in U$. 
Hence, noting that $\mathbb{R}^m\times \mathbb{R}_{++}\subset \mathrm{dom}(S)$, 
by Lemmas \ref{sfj} and \ref{L516}, $S$ is outer Lipschitzian on $U$ with modulus $\max\limits_{J_1,J_2\subset  [n],~J_1\cap J_2=\emptyset} \kappa_{J_1 J_2} L_H$.
\end{proof}

Using a notable result on Lipschitz continuity from Robinson \cite{Robinson07}, the next step to prove the Lipschitz continuity of $S$ is to demonstrate its inner semicontinuity. 

\begin{lemma}(\cite[Theorem 1.5]{Robinson07})\label{L516a}
Let $\Phi:\mathbb{R}^p \rightrightarrows \mathbb{R}^q$ be a multifunction having closed values, and let $D$ be a convex subset of $\mathrm{dom}(\Phi)$, and let $\mu$ be a nonnegative real number. Then, $\Phi$ is Lipschitz continuous relative to $D$ with modulus $\mu$ if and only if $\Phi$ is outer Lipschitz continuous relative to $D$ with modulus $\mu$ and is inner semicontinuous relative to $D$.
\end{lemma}


To prove the inner semicontinuity of $S$, we show an important consistency property of the decomposition \eqref{510f}. 

\begin{proposition}\label{lemma-S-lip-inner}
Suppose that $H(b,\lambda)$ is single-valued on an open set $U\subset \mathbb{R}^m\times \mathbb{R}_{++}$, 
and $h$ is continuously differentiable at $Ax-b$ for all $(b,\lambda)\in U$ with $Ax\in H(b,\lambda)$. 
Then, 
\begin{equation}\label{0429a}
  S( b, \lambda)=S_{J_1(b,\lambda)J_2(b,\lambda)}( b, \lambda)
  ~~\forall (b,\lambda)\in U,
\end{equation}
where
\begin{equation}\label{index-J1-J2}
\begin{aligned}
    J_1(b,\lambda):=&\{i\in [n]: \frac{1}{\lambda}A_i^T\nabla h(Ax-b)=-1~\mbox{with}~Ax\in H(b,\lambda)\},\\
    J_2(b,\lambda):=&\{i\in [n]: \frac{1}{\lambda}A_i^T\nabla h(Ax-b)=1~\mbox{with}~Ax\in H(b,\lambda)\}.
\end{aligned}   
\end{equation}
Moreover, if $\lim\limits_{k\rightarrow\infty}(b^k,\lambda^k) = (\widehat b,\widehat\lambda)\in U$ and $S( b^k,\lambda^k)= S_{\tilde J_1\tilde J_2}( b^k,\lambda^k)$ for all $k$. Then
\begin{equation}\label{s=}
    S(\widehat b,\widehat\lambda)=S_{\tilde J_1\tilde J_2}(\widehat b,\widehat\lambda).
\end{equation}
\end{proposition}

\begin{proof}
First, for each $(b,\lambda)\in U$, since $H(b,\lambda)$ is a singleton, the indexes $J_1(b,\lambda)$ and $J_2(b,\lambda)$ are well defined. 
Moreover, for any $x\in S(b,\lambda)$, since $h$ is continuously differentiable at $Ax-b$, by \eqref{e101}, we have 
$
    -\frac{1}{\lambda}A^T\nabla h(Ax-b)\in \partial \|\cdot\|_1(x).
$
Thus, $x_i\ge 0$ for $i\in J_1(b,\lambda)$, $ x_j\leq 0$ for $j\in J_2(b,\lambda)$, 
$x_k=0$ for $k\in (J_1(b,\lambda)\cup J_2(b,\lambda))^C$. 
Combining this with the definition of $J_1(b,\lambda),~J_2(b,\lambda)$, we have 
\begin{equation*}
x\in P_{J_1(b,\lambda)J_2(b,\lambda)}
~~\mathrm{and}~~
    -\frac{1}{\lambda}A^T\nabla h(Ax-b)\in  E_{J_1(b,\lambda)J_2(b,\lambda)}. 
\end{equation*}
This implies that $S( b, \lambda)\subset S_{J_1(b,\lambda)J_2(b,\lambda)}( b, \lambda)$. Then \eqref{0429a} follows from Lemma \ref{sfj}.

Second, by \eqref{0429a}, $S(\widehat b,\widehat \lambda)=S_{\widehat J_1\widehat J_2}(\widehat b,\widehat \lambda)$, where $\widehat J_j:=J_j(\widehat b,\widehat \lambda)$ for $j=1,2$.  
Take $x^k\in S(b^k,\lambda^k)=S_{\tilde J_1\tilde J_2}( b^k,\lambda^k)$. 
Since $\widehat\lambda > 0$, by Lemma \ref{L625a}, the sequence $\{x^k\}_{k=1}^\infty$ is bounded. Without loss of generality, we can assume that $x^k \rightarrow \tilde{x}$ as $k\rightarrow\infty$ (by taking a subsequence if necessary). By Lemma \ref{L625a}, $\tilde x\in S(\widehat b,\widehat \lambda)=S_{\widehat J_1\widehat J_2}(\widehat b,\widehat \lambda)$.  
On the other hand, since $x^k\in S_{\tilde J_1\tilde J_2}( b^k,\lambda^k)$, we also have $\tilde x\in S_{\tilde J_1\tilde J_2}( \widehat b,\widehat \lambda)$ due to \eqref{510g} and the continuity of $\nabla h$ at $A\tilde x-\widehat b$. Hence, $\tilde J_1\subset J_1(\widehat b,\widehat\lambda)=\widehat J_1$ and $\tilde J_2\subset J_2(\widehat b,\widehat \lambda)=\widehat J_2$ by the definition of $E_{\tilde J_1\tilde J_2}$ in \eqref{ePE} and $S_{\tilde J_1\tilde J_2}$ in \eqref{510g}. It follows that
\begin{equation}\label{0429c}
\{x^k\}_{k=1}^\infty\subset  P_{\tilde J_1\tilde J_2}\subset P_{\widehat J_1\widehat J_2}~\text{and}~\tilde x\in P_{\tilde J_1\tilde J_2}\subset P_{\widehat J_1\widehat J_2}.
\end{equation}
For any $\widehat{x}\in S(\widehat b,\widehat \lambda)$, 
since $H(\widehat b,\widehat \lambda)$ is single-valued and $\tilde x\in S(\widehat b,\widehat \lambda)$, $A\widehat{x}=A\tilde{x}$. Then $h(A\widehat{x}-b)=h(A\tilde{x}-b)$. Thus $\|\widehat{x}\|_1=\|\tilde{x}\|_1$. 
Since $\widehat{x}\in  P_{\widehat J_1\widehat J_2}$, by \eqref{0429c}, we have
\begin{equation}\label{0429b}
  \sum\limits_{i\in \widehat J_1}\widehat{x}_i-\sum\limits_{j\in \widehat J_2}\widehat{x}_j=\sum\limits_{i\in \widehat J_1}\tilde{x}_i-\sum\limits_{j\in \widehat J_2}\tilde{x}_j.
\end{equation}
Since $x^k\rightarrow \tilde x$, for all sufficiently large $k$, 
$\mathrm{sgn}\left(x^k-\frac{\tilde{x}}{2}\right)_i=\mathrm{sgn}\left(\tilde{x}\right)_i$ for all $i\in \mathrm{supp}(\tilde{x})$. 
Then, by \eqref{0429c}, we have $\left(x^k-\frac{\tilde{x}}{2}\right)\in P_{\widehat J_1\widehat J_2}$. Thus $x^k+\frac{\widehat{x}-\tilde{x}}{2}\in P_{\widehat J_1\widehat J_2}$. Then 
\begin{equation*}
  \aligned
  \left\|x^k+\frac{\widehat{x}-\tilde{x}}{2}\right\|_1&=\sum\limits_{i\in \widehat J_1}\left(x_i^k+\frac{\widehat{x}_i-\tilde{x}_i}{2}\right)-\sum\limits_{j\in \widehat J_2}\left(x_j^k+\frac{\widehat{x}_j-\tilde{x}_j}{2}\right)~\\
  &=\sum\limits_{i\in \widehat J_1}x_i^k-\sum\limits_{j\in \widehat J_2}{x}_j^k+\sum\limits_{i\in \widehat J_1}\frac{(\widehat{x}_i-\tilde{x}_i)}{2}-\sum\limits_{j\in \widehat J_2}\frac{(\widehat{x}_j-\tilde{x}_j)}{2}\\
  &=\sum\limits_{i\in \widehat J_1}x_i^k-\sum\limits_{j\in \widehat J_2}{x}_j^k
  =\|x^k\|_1 ~~~~~~~~(\text{by}~ \eqref{0429b}~\mathrm{and}~\eqref{0429c}) .
  \endaligned
\end{equation*}
Noting that $A\left(x^k+\frac{\widehat{x}-\tilde{x}}{2}\right)=A(x^k)$ (thanks to the fact $A\widehat{x}=A\tilde{x}$),  this implies that
$$h\left(A\left(x^k+\frac{\widehat{x}-\tilde{x}}{2}\right)-b^k\right)+\lambda^k\left\|x^k+\frac{\widehat{x}-\tilde{x}}{2}\right\|_1=h(Ax^k-b^k)+\lambda^k\|x^k\|_1.$$
Since $x^k\in S(b^k,\lambda^k)$, 
it follows that $x^k+\frac{\widehat{x}-\tilde{x}}{2}\in S(b^k,\lambda^k)\subset P_{\tilde J_1\tilde J_2}$. 
Noting that $x^k,~\tilde x\in P_{\tilde J_1\tilde J_2}$, it follows that $\widehat{x}_i=2(x^k+\frac{\widehat{x}-\tilde{x}}{2} - x^k +\frac{\tilde{x}}{2})_i=0$ for all $i\in (\tilde J_1\cup\tilde J_2)^C$. Noting that $\widehat{x}\in P_{\widehat J_1\widehat J_2}$, $\tilde J_1\subset \widehat J_1$ and $\tilde J_2\subset \widehat J_2$, we have $\widehat{x}\in P_{\tilde J_1\tilde J_2}$ and $E_{\widehat J_1\widehat J_2}\subset E_{\tilde J_1\tilde J_2}$. Hence, $\widehat{x}\in S_{\tilde J_1\tilde J_2}(\widehat b,\widehat\lambda)$ by \eqref{510g} and using the fact $\widehat{x}\in S_{\widehat J_1\widehat J_2}(\widehat b,\widehat\lambda)$. Thus $S(\widehat b,\widehat\lambda)\subset S_{\tilde J_1\tilde J_2}(\widehat b,\widehat\lambda)$. 
Therefore, the equality \eqref{s=} follows from Lemma \ref{sfj}.
\end{proof}

\begin{theorem}\label{T01}
(i) Suppose that $H$ is single-valued and Lipschitz continuous on a convex open set $U\subset  \mathbb{R}^m\times \mathbb{R}_{++}$, 
and $h$ is continuously differentiable at $Ax-b$ for all $(b,\lambda)\in U$ with $Ax\in H(b,\lambda)$. 
Then $S$ is Lipschitz continuous relative to $U$.

(ii) Let $h:\mathbb{R}^m\rightarrow \mathbb{R}_+$ be a continuous convex function. Given $(\bar b,\bar\lambda)\in \mathbb R^m\times \mathbb R_{++}$, suppose that Condition \ref{conditionH1-0} holds. 
Then the solution mapping $S$ of \eqref{Psc} is locally Lipschitz continuous at $(\bar b, \bar \lambda)$. 
\end{theorem}

\begin{proof}
(i) By Lemma \ref{L516a} and Corollary \ref{S-outer Lipschitzian}, it suffices to prove the inner semicontinuity of $S$ relative to $U$.  
Suppose not, there exists $(\widehat b,\widehat\lambda)\in U$ such that $S$ is not inner semicontinuous at $(\widehat b,\widehat\lambda)$ relative to $U$, i.e., 
$
  S(\widehat b,\widehat\lambda)\not\subset \liminf\limits_{(b,\lambda)\rightarrow (\widehat b,\widehat\lambda)}S(b,\lambda).
$ It follows that there exists $\widehat{x}\in S(\widehat b,\widehat \lambda)$ such that $\widehat{x}\notin \liminf\limits_{(b,\lambda)\rightarrow (\widehat b,\widehat\lambda)}S(b,\lambda)$. Hence, there is a positive constant $\alpha$ and a sequence $(b^k,\lambda^k)\rightarrow (\widehat b,\widehat\lambda)$ such that 
\begin{equation}\label{e626a}
    d(\widehat{x},S(b^k,\lambda^k))\ge \alpha >0
    ~~\forall k\in \mathbb N.
\end{equation}
Since $U$ is open, for sufficiently large $k$, by \eqref{0429a} in Proposition \ref{lemma-S-lip-inner}, we have $S(b^k,\lambda^k)=S_{J_1(b^k,\lambda^k) J_2(b^k,\lambda^k)}(b^k,\lambda^k)$.
Note that the set $[n]$ is finite, so by taking a subsequence if necessary, we can assume without loss of generality that 
$J_1(b^k,\lambda^k)=\tilde{J}_1\subset  [n]$ and $J_2(b^k,\lambda^k)=\tilde{J}_2\subset  [n]$ for all $k$. 
Hence, $S(b^k,\lambda^k)=S_{\tilde J_1\tilde J_2}(b^k,\lambda^k)$. 
Proposition \ref{lemma-S-lip-inner} shows that 
$S(\widehat b,\widehat\lambda)=S_{\tilde J_1\tilde J_2}(\widehat b,\widehat\lambda)$. Thus $\widehat{x}\in S_{\tilde J_1\tilde J_2}(\widehat b,\widehat\lambda)$. 
Since $S_{\tilde J_1\tilde J_2}$ is Lipschitzian on $\mathrm{dom}(S_{\tilde J_1\tilde J_2})\cap  U$ by Lemma \ref{lemma-S-lip}, it follows that $$ 
d(\widehat{x},S(b^k,\lambda^k))\leq d_H(S_{\tilde J_1\tilde J_2}(\widehat b,\widehat \lambda),S(b^k,\lambda^k))=d_H(S_{\tilde J_1\tilde J_2}(\widehat b,\widehat \lambda),S_{\tilde J_1\tilde J_2}(b^k,\lambda^k))\rightarrow 0.
$$
This contradicts \eqref{e626a}. Hence, $S$ is inner semicontinuous relative to $U$. 
Therefore, $S$ is Lipschitz continuous relative to $U$. 

(ii) Suppose Condition \ref{conditionH1-0} holds, by Theorem \ref{aTi1}, $H$ is single-valued and Lipschitz continuous around $(\bar b,\bar \lambda)$. By (i), $S$ is Lipschitz continuous at $(\bar b, \bar \lambda)$.
\end{proof}

\begin{remark}\label{Lconstant}
   By Lemma \ref{L516a} and Corollary \ref{S-outer Lipschitzian}, we know that the Lipschitz constant of $S$ is $\max_{J_1,J_2\subset  [n], J_1\cap J_2=\emptyset}\kappa_{J_1J_2} L_H$, 
   where $\kappa_{J_1J_2}$ and $L_H$ are the Lipschitz constants of $F_{J_1J_2}$ and $H$, respectively. Note that $L_H$ can be calculated from \eqref{lipcon} using \eqref{e72}.
\end{remark}


\subsection{Uniqueness, Lipschitz and smoothness properties of solution mapping}

In this section, we mainly establish a sufficient condition for the local single-valued and Lipschitz continuity of solution mapping for the LASSO-type problem \eqref{Psc}.
To this end, we first introduce a sufficient condition for solution uniqueness.

\begin{proposition}\label{Pro1}
Let $\bar x$ be a solution of \eqref{Psc}.
Suppose that $h$ is continuous differentiable at $A\bar x-b$, and $H(b,\lambda)$ is a singleton, namely, $\{Ax:x\in S(b,\lambda)\}=\{A\bar x\}$. Assume also that $\ker (A_{J_{h}})=\{0\}$, or equivalently,  $A_{J_{h}}$ has full column rank, where $J_{h}:=\left\{i\in [n]:|A_i^T\nabla h(A\bar x- b)|=\lambda\right\}$. 
Then  $\bar x$ is the unique solution of \eqref{Psc}.
\end{proposition}

\begin{proof}
For any solution $x\in S(b,\lambda)$, since $Ax=A\bar x$, by the definition of $J_h$ and the optimality condition \eqref{e101} of $\bar x\in S(b,\lambda)$, we have
$$A_i^T\nabla h(A x-b)=A_i^T\nabla h(A\bar x-b)\in (-\lambda,\lambda)~~\forall i\in J_h^C.$$
Hence, by \eqref{e101}, $x_{J_h^C}=0=\bar x_{J_h^C}$. 
Thus $A_{J_h}(x_{J_h})=Ax=A\bar x=A_{J_h}(\bar x_{J_h})$. Since $\ker (A_{J_{h}})=\{0\}$, we get $x_{J_h}=\bar x_{J_h}$. Therefore, $x=\bar x$ and the proof is complete. 
\end{proof}

Next we establish the Lipschitz and smoothness properties of $S$ under proper conditions. 

\begin{theorem}\label{T318-h}
Let $h:\mathbb{R}^m\rightarrow \mathbb{R}_{+}$ be a continuous convex function.  \\
(i) If Condition \ref{conditionH2-0} holds at $(\bar b, \bar \lambda)\in\mathbb{R}^m\times\mathbb{R}_{++}$, then~the solution mapping $S$ of \eqref{Psc} is locally single-valued and Lipschitz around $(\bar b,\bar \lambda)$.\\
(ii) If Condition \ref{conditionH2-0} holds at $(\bar b, \bar \lambda)\in\mathbb{R}^m\times\mathbb{R}_{++}$, 
and $J_h=\mathrm{supp}(\bar x)$ for $\bar x\in S(\bar b,\bar \lambda)$,
then $S$ is continuously differentiable at $(\bar b,\bar \lambda)$.
 \end{theorem}
 
 \begin{proof}
(i) Since Condition \ref{conditionH2-0} is equivalent to Condition \ref{conditionH1-0} plus $A_{J_{h}}$ has full column rank. 
By Theorem \ref{aTi1}, it remains to prove that $S$ is single-valued around $(\bar b,\bar \lambda)$. In the proof of Theorem \ref{aTi1}, we have proved that the mapping $H (b,\lambda)$ is single-valued and Lipschitz continuous on $B(\bar b,\delta_4)\times B(\bar \lambda ,\delta_4)$. 
By \eqref{510e}, 
\begin{equation*}
|A_i^T\nabla h(A x- b)|<\lambda~~\forall i\in  J_h^C,~\forall x\in S(b,\lambda), 
~\forall (b,\lambda)\in B(\bar b,\delta_4)\times B(\bar\lambda,\delta_4) .
\end{equation*}
This implies that
\begin{equation}\label{index-conlusion}
J_{h(b,\lambda)}:=\left\{i\in [n]:|A_i^T\nabla h(A x- b)|=\lambda~\text{for}~x\in S(b,\lambda)\right\}\subset J_h. 
\end{equation}
Since $A_{J_{h}}$ has full column rank, $A_{J(b,\lambda)}$ does too. 
By Proposition \ref{Pro1}, the solution set $S(b,\lambda)$ is a singleton. Hence, $S$ is single-valued on $B(\bar b,\delta_4)\times B(\bar \lambda ,\delta_4)$.

 (ii) First, assume that $J_h\neq \emptyset$. 
Since $\bar \lambda>0$ and $h$ is twice continuously differentiable in a neighborhood of $A\bar x-\bar b$, 
there is a positive constant $\alpha_1$ such that 
\begin{equation*}
f: 
(x_{J_h},b,\lambda)
\mapsto\frac{1}{\lambda} A_{J_h}^T\nabla h(A_{J_h} x_{J_h}-b) +d 
\end{equation*}
is well defined for $(x_{J_h},b,\lambda)\in  B(\bar x_{J_h},\alpha_1)\times B(\bar b,\alpha_1)\times B(\bar \lambda,\alpha_1) \subset \mathbb{R}^{|J_h|}\times \mathbb{R}^m\times\mathbb{R}_{++}$ (noting that $A\bar x= A_{J_h}\bar x_{J_h}$). Here $d\in \mathbb{R}^{|J_h|}$ with $d_i=1$ for $i\in \{i\in  [n]:A_i^T\nabla h(A\bar x-\bar b)=-\bar\lambda \}$ and $d_j=-1$ for $j\in \{j\in  [n]:A_j^T\nabla h(A\bar x-\bar b)=\bar\lambda \}$.
Since $\bar x_{J_h^C}=0$, $f(\bar x_{J_h},\bar b,\bar \lambda)=0$ and $D_{x_{J_h}}f(\bar x_{J_h},\bar b,\bar \lambda)=\frac{1}{\bar \lambda} A_{J_h}^T\nabla^2 h(A\bar x-b) A_{J_h}$. 
Since $h$ is convex, by the fact $\ker(\nabla^2h(A\bar x- \bar b))\cap \mathrm{rge}(A_{J_h})=\{0\}$, and $\ker(A_{J_h})=\{0\}$, we have 
\begin{equation*}
    \aligned
\ker(D_{x_{J_h}}f(\bar x_{J_h},\bar b,\bar \lambda))=&\{x_{J_h}\in \mathbb{R}^{|J_h|}:A_{J_h}^T\nabla^2 h(A\bar x-b) A_{J_h} x_{J_h}=0\}\\
=&\{x_{J_h}\in \mathbb{R}^{|J_h|} :x_{J_h}^TA_{J_h}^T\nabla^2 h(A\bar x-b) A_{J_h}x_{J_h}=0\}\\
=&\{x_{J_h}\in \mathbb{R}^{|J_h|} :\nabla^2 h(A\bar x-b) A_{J_h}x_{J_h}=0\}\\
=&\{x_{J_h}\in \mathbb{R}^{|J_h|} : A_{J_h}x_{J_h}=0\}=\{0\} ,
    \endaligned
\end{equation*}
where the third equality comes from the argument of the proof of \eqref{520f}. Hence, $D_{x_{J_h}} f(\bar x_{J_h},\bar b,\bar \lambda)$ is nonsingular. By the classical implicit function theorem, there exists $\alpha_2\in (0, \alpha_1)$ such that for each $(b,\lambda)\in B(\bar b,\alpha_2)\times B(\bar\lambda,\alpha_2)$, there is a unique point $\widehat S(b,\lambda)\in B(\bar x_{J_h},\alpha_2)$ satisfying
$
  f(\widehat S(b,\lambda),b,\lambda)=0.
$
 Moreover, $\widehat S$ is continuously differentiable in $B(\bar b,\alpha_2)\times B(\bar\lambda,\alpha_2)$ due to the continuous differentiability of $f$. 

 On the other hand, by the conclusion of (i), $S(b,\lambda)$ is single-valued and locally Lipschitz on $B(\bar b,\delta_4)\times B(\bar\lambda,\delta_4)$. 
 Moreover, by \eqref{index-conlusion} and \eqref{e101}, we have $S(b,\lambda)_{J_h^C}=0$. 
 (Here we do not distinguish between a singleton set $S(b,\lambda)$ and its unique element for simplicity.)  
 Note that $S(\bar b, \bar \lambda )=\{\bar x\}$. 
 By the Lipschitz continuity of $S$ and the assumption $\mathrm{supp}(\bar x)=J_h$, there exists $\alpha_3\in (0,\min\{\alpha_2,\delta_4\})$ such that $\mathrm{sgn}(S(b,\lambda)_{J_h})=\mathrm{sgn}(\bar x_{J_h})$ 
 and $S(b, \lambda )\in B(\bar x_{J_h},\alpha_2)$
 for all $(b,\lambda)\in B(\bar b, \alpha_3)\times B(\bar\lambda, \alpha_3)$. Hence, for such $(b,\lambda)$,  
$f(S(b,\lambda)_{J_h},b,\lambda)=0$ by \eqref{e101}. 
Therefore, by the uniqueness, we have $S(b,\lambda)_{J_h}\equiv \widehat S(b,\lambda)$. Recall that $S(b,\lambda)_{J_h^C}=0$. It is proven that $S(b,\lambda)$ is continuously differentiable in $B(\bar b, \alpha_3)\times B(\bar\lambda, \alpha_3)$.

Second, we consider the case where $J_h= \emptyset$. 
  In this case, by \eqref{index-conlusion} and \eqref{e101}, for all $(b,\lambda)\in B(\bar b,\delta_4)\times B(\bar \lambda ,\delta_4)$ and all $x\in S(b,\lambda)$, we have $x=0$. Namely, $S(b,\lambda)\equiv 0$ in $(\bar b,\delta_4)\times B(\bar \lambda ,\delta_4)$. Hence, it is continuously differentiable.
 \end{proof}

\section{The LASSO problem}
\label{section-lasso}
\setcounter{equation}{0}

In this section, we focus on the LASSO problem and explore the benefits of properly utilizing its structure.

\subsection{Multi-valued case for LASSO}

This section establishes the global Lipschitz continuity of solution mapping for the LASSO. It leverages the polyhedral property of  $S$ for the LASSO and a general result of Robinson \cite{Robinson07} on polyhedral multifunctions, given below. 

\begin{lemma}(\cite[Corollary 2.2]{Robinson07})\label{Ls}
Let $F$ be a polyhedral multifunction from $\mathbb{R}^n$ to $\mathbb{R}^m$, let $\mu$ be its outer Lipschitzizan modulus,  and let $D$ be a convex subset of $\mathrm{dom}(F)$ on which $F$ is single-valued. 
Then $F$ is Lipschitzian on $D$ with modulus $\mu$.
\end{lemma}

\begin{proposition}\label{lemma-S-polyhedral}
For the LASSO problem \eqref{P1}, the solution mapping $S$ is a polyhedral multifunction.
\end{proposition}

\begin{proof}
    By Lemma \ref{sfj}, $\mathrm{gph}(S)=\bigcup\limits_{J_1,J_2\subset  [n],~J_1\cap J_2=\emptyset} \mathrm{gph}(S_{J_1 J_2})$, where 
    $$ S_{J_1 J_2}(b,\lambda):=\left\{x\in P_{J_1 J_2}: 0\in  {A^T(Ax-b)}+ {\lambda} E_{J_1 J_2}\right\}.$$
    Then we have
    \begin{equation*}
    \aligned
        \mathrm{gph}(S_{J_1 J_2})=&\{(b,\lambda,x)\in \mathbb R^m\times \mathbb R_+\times P_{J_1 J_2}:- {A^T(Ax-b)}\in{\lambda} E_{J_1 J_2}\}\\
        =&\{(b,\lambda,x)\in \mathbb R^m\times \mathbb R_+\times \mathbb{R}^n : A_i^TAx-A_i^Tb+ \lambda=0, ~x_i\geq 0~~\forall i\in J_1,\\
       &~~~~~~~~~~~~~~~~~~~~ ~~~~~~~~~~~~~~ A_j^TAx-A_j^Tb- \lambda=0, ~x_j\leq 0 ~~\forall j\in J_2,\\
       &~~~~~~~~~~~~~~~~~~~~  -\lambda\leq A_k^TAx-A_k^Tb\leq \lambda,
       ~x_k= 0~~\forall k\in (J_1\cup J_2)^C\}.
        \endaligned
    \end{equation*}
    Hence $S$ is a polyhedral multifunction.
\end{proof}


\begin{theorem}\label{Ti1}
For any $A\in \mathbb{R}^{m\times n}$, the solution mapping  $S$ of the LASSO problem \eqref{P1}
is globally (Hausdorff) Lipschitz continuous on $\mathbb{R}^m\times \mathbb{R}_{++}$.
\end{theorem}

\begin{proof}
By Proposition \ref{lemma-S-polyhedral}, $S$ is a polyhedral multifunction. Hence, $H$ defined in \eqref{1521a-0} is also a polyhedral multifunction. By \cite[Theorem 208]{Shapiro00}, there exists $L_H>0$ such that $H$ is outer Lipschitzian at every $(b,\lambda)\in \mathbb{R}^m\times \mathbb{R}_{++}$ with modulus $L_H$.
Since $H$ is single-valued on $\mathbb{R}^m\times \mathbb{R}_{++}$ by Proposition \ref{L514a}, it follows from Lemma \ref{Ls} that $H$ is Lipschitz on $\mathbb{R}^m\times \mathbb{R}_{++}$. Using Theorem \ref{T01} (i), the solution mapping $S$ is globally Lipschitz on $\mathbb{R}^m\times \mathbb{R}_{++}$.
\end{proof}

\subsection{Single-valued case for LASSO}

In this section, we fully characterize the local single-valued and Lipschitz continuity of solution mappings for the LASSO using Tibshirani’s sufficient condition. Indeed, this can be done for the general problem \eqref{Psc} with a strongly convex loss function.

\begin{theorem}\label{T318}
Let $h:\mathbb{R}^m\rightarrow \mathbb{R}_{+}$ be a strongly convex, twice continuously differentiable function. 
Then, the solution mapping $S$ of \eqref{Psc} is locally single-valued and Lipschitz around $(\bar b,\bar \lambda)$
if and only if $A_{J_h}$ has full column rank. Here, $J_{h}:=\left\{i\in [n]:|A_i^T\nabla h(A\bar x-\bar b)|=\bar\lambda\right\}$ is defined in \eqref{J-LASSO}, where $\bar x\in S(\bar b,\bar \lambda)$. 
 \end{theorem}
\begin{proof}
The strong convexity of $h$ implies that $ \ker (\nabla^2h(A\bar x-\bar b))=\{0\}$. Thus, Condition \ref{conditionH2-0} is valid. Hence, Theorem \ref{T318-h} (i) establishes the sufficiency. 
It remains to prove the necessity. Take $\bar x\in S(\bar b,\bar \lambda)$ and denote
$$J_1:=\left\{i\in  [n]:A_i^T\nabla h(A\bar x-\bar b)=-\bar\lambda\right\}~
\text{and}~
J_2:=\left\{i\in  [n]:A_i^T\nabla h(A\bar x-\bar b)=\bar\lambda\right\}.$$
Define $a\in \mathbb{R}^n$ as \eqref{ea}. Suppose that $S$ is locally single-valued and Lipschitz around $(\bar b,\bar \lambda)$. Then there exists $t>0$ such that $S$ is single-valued at $(\bar b+t Aa,\bar\lambda)$. By Lemma \ref{L613a}, $\mathrm{supp}(\bar x+t a)=J_1\cup J_2=J_h$ and $\bar x+t a\in S(\bar b+t Aa,\bar \lambda)$. Since $S$ is single-valued at $(\bar b+t Aa,\bar\lambda)$, the solution
$\bar x+t a$ of \eqref{P1} with $(b,\lambda)=(\bar b+t Aa,\bar\lambda)$ is unique. Note that  $\mathrm{supp}(\bar x+t a)=J_h$, by \cite[Theorem 2.1]{ZY}, $A_{J_h}$ has full column rank.
 \end{proof}


As a corollary of Theorem \ref{T318}, the following theorem fully characterizes the local single-valued and Lipschitz continuity of solution mapping for the LASSO.

\begin{theorem}\label{T318-lasso}
The solution mapping $S$ of the LASSO problem \eqref{P1} is locally single-valued and Lipschitz around $(\bar b,\bar \lambda)$
if and only if Tibshirani's sufficient condition holds, \textit{that is}, $A_{J_{LA}}$ has full column rank, or equivalently, $\ker (A_{J_{LA}})=\{0\}$, where \begin{equation}\label{J-LASSO-0}
    J_{LA}=J_{LA}(\bar x):=\{i\in [n]:|A_i^T( \bar b-A\bar x)|= \bar \lambda\}.
\end{equation}
 \end{theorem}
\begin{remark}
It is noteworthy that Nghia \cite[Theorem 3.12]{nghia2024geometric} recently demonstrated that Tibshirani’s sufficient condition is also necessary, employing the second subderivative characterization of strong minima from \cite[Theorem 13.24]{rockwets98}. In contrast, our approach relies solely on the necessary result of solution uniqueness for the LASSO, as shown by Zhang \textit{et al.} \cite[Theorem 2.1]{ZY}.
\end{remark}
 
\section{The SR-LASSO problem}
\label{section-sr-lasso}
\setcounter{equation}{0}

Compared to LASSO, the Lipschitz behavior for SR-LASSO is more complex. In this section, we will fully characterize the Lipschitz property of solution mapping for the SR-LASSO problem \eqref{Pi1} in both multi-valued and single-valued cases.

\subsection{Multi-valued case for SR-LASSO}

The Lipschitz continuity of solution mapping for the SR-LASSO requires certain conditions since $S_{SR}$ in Example \ref{fla} is not Lipschitz continuous at $(b, 1)$.   
This section is devoted to establishing a characterization of Lipschitz continuity of solution mapping for the SR-LASSO. 

We first show that Condition \ref{conditionH1-0} is equivalent to Condition \ref{A61} for the SR-LASSO.

\begin{proposition}\label{LCsr}
    For $h(z)=\|z\|$, Condition \ref{conditionH1-0} is equivalent to Condition \ref{A61}.
\end{proposition}

\begin{proof}
Clearly, $h(z)=\|z\|$ being twice continuously differentiable around $A\bar x - \bar b$ is equivalent to $A\bar x \neq \bar b$. 
For $A\bar x\neq \bar b$, we have 
 $\nabla h(A\bar x-\bar b)=\frac{A\bar x-\bar b}{\|A\bar x-\bar b\|},$
 and then $J_{SR}(\bar x)=J_{h}(\bar x)$. 
 The remaining part of the proof is to show that 
 \begin{equation}\label{e619}
 \bar b\notin \mathrm{rge}(A_{J_{SR}})
 \iff 
     \ker(\nabla^2h(A\bar x- \bar b))\cap \mathrm{rge}(A_{J_h})=\{0\},~\mathrm{for}~A\bar x\neq \bar b.
 \end{equation}
In fact, suppose that $\bar b\notin \mathrm{rge}(A_{J_{SR}})$. Since $A\bar x \in \mathrm{rge}(A_{J_{SR}})$ by Lemma \ref{L613a} (i), we have $\frac{(A\bar x-\bar b)}{\|A\bar x-\bar b\|}\notin \mathrm{rge}(A_{J_{SR}})$. Thus, for any $0\neq u\in \mathrm{rge}(A_{J_{SR}}) $, we have 
\begin{equation}\label{e630-a}
    \frac{u}{\|u\|}\neq\frac{(A\bar x-\bar b)}{\|A\bar x-\bar b\|}~\text{and}~\left|\left\langle\frac{A\bar x- \bar b}{\|A\bar x-\bar b\|},u\right\rangle\right|<\left\|\frac{A\bar x- \bar b}{\|A\bar x- \bar b\|}\right\|\|u\|=\|u\|.
\end{equation}
Note that 
\begin{equation}\label{s2}
    \nabla^2 h(A\bar x-\bar b)=\frac{1}{\|A\bar x-\bar b\|}\left(\mathbb{I}-\frac{A\bar x-\bar b}{\|A\bar x-\bar b\|}\left(\frac{A\bar x-\bar b}{\|A\bar x-\bar b\|}\right)^T\right),
\end{equation}
where $\mathbb{I}$ is the identity matrix. 
It follows that
\begin{equation*}
  \aligned
  \left\langle \left[\nabla^2 h(A\bar x-\bar b)\right]u, u\right\rangle 
  =\frac{1}{\|A\bar x- \bar b\|}\left(\|u\|^2-\left|\left\langle\frac{A\bar x- \bar b}{\|A\bar x- \bar b\|},u\right\rangle\right|^2\right) >0
  ~~~(\text{by}~\eqref{e630-a}) .
  \endaligned
\end{equation*}
This implies that $\left[\nabla^2 h(A\bar x-\bar b)\right]u\neq0$. Hence the ``$\Longrightarrow$" of \eqref{e619} is proved. 

Conversely, suppose not, then $\bar b\in \mathrm{rge}(A_{J_{SR}})$. Thus $A\bar x-\bar b\in \mathrm{rge}(A_{J_{SR}})$ by Lemma \ref{L613a} (i). 
On the other hand, by \eqref{s2}, we have 
\begin{equation*}
  \aligned
  \left[\nabla^2 h(A\bar x-\bar b)\right](A\bar x-\bar b) 
  = \frac{1}{\|A\bar x- \bar b\|}\left((A\bar x-\bar b)-(A\bar x-\bar b)\right)=0.
  \endaligned
\end{equation*}
This implies that $A\bar x-\bar b\in \ker(\nabla^2h(A\bar x- \bar b))$. 
Hence $0\neq A\bar x-\bar b\in \ker(\nabla^2h(A\bar x- \bar b))\cap \mathrm{rge}(A_{J_h})$. 
This contradicts the condition on the right-hand side of \eqref{e619}. 
\end{proof}

With the help of Proposition \ref{LCsr} and Theorem \ref{T01}, we can establish sufficient and necessary condition for the local Lipschitz continuity of solution mapping for the SR-LASSO problem \eqref{Pi1}.

\begin{theorem}\label{ca3}
Let $(\bar b, \bar \lambda)\in\mathbb{R}^m\times\mathbb{R}_{++}$. Suppose that $\{Ax:x\in S(\bar b,\bar\lambda)\}\neq \{\bar b\}$. Then the solution mapping $S$ of the SR-LASSO
is locally Lipschitzian at $(\bar b,\bar \lambda)$ if and only if Condition \ref{A61} holds.
 \end{theorem}
 
 \begin{proof}
 \textit{Sufficiency.} 
 If Condition \ref{A61} holds, by Proposition \ref{LCsr}, Condition \ref{conditionH1-0} holds for $h(z)=\| z \|$.  
 Using Theorem \ref{T01} (ii), $S$ is locally Lipschitzian at $(\bar b,\bar \lambda)$. 

\textit{Necessity.} 
Suppose that $S$ is locally Lipschitz at $(\bar b,\bar \lambda)$, but Condition \ref{A61} is not fulfilled. Then $\bar b\in \mathrm{rge}(A_{J_{SR}})$ for $\bar x\in S(\bar b,\bar \lambda)$ satisfying $A\bar x \neq \bar b$ (noting that $\{Ax:x\in S(\bar b,\bar\lambda)\}\neq \{\bar b\}$). 
Denote
$$J_1:=\left\{i\in  [n]:A_i^T\frac{A\bar x -\bar b}{\bar\lambda\|A\bar x-\bar b\|}= -1\right\}~\text{and}~
J_2:=\left\{i\in  [n]:A_i^T\frac{A\bar x -\bar b}{\bar\lambda\|A\bar x-\bar b\|}= 1\right\},$$
and define $a\in \mathbb R^n$ as \eqref{ea}. 
By the assumption that $S$ is locally Lipschitz at $(\bar b,\bar \lambda)$, there exists $t_0>0$ sufficiently small such that $S$ is locally Lipschitz at $(\bar b+t_0Aa,\bar \lambda)$. 
By Lemma \ref{L613a} (ii), $\mathrm{sgn}(\bar x+t_0a)=\mathrm{sgn}(a)$ and  $\bar x+t_0a\in S(\bar b+t_0Aa,\bar \lambda).$ 
Since $S$ is locally Lipschitz at $(\bar b+t_0Aa,\bar \lambda)$, we have 
 $\lim\limits_{\lambda\rightarrow \bar\lambda}d(\bar x+t_0a,S(\bar b+t_0Aa,\lambda))=0$. 
 Hence, there exists $x^k\in S(\bar b+t_0Aa,\bar\lambda+1/ k)$ such that $x^k\rightarrow \bar x+t_0a$ (as $k\rightarrow \infty$). 
 Recalling that $\mathrm{sgn}(\bar x+t_0a)=\mathrm{sgn}(a)$, we have $(x^k)_i>0$ for $i\in J_1$ and $(x^k)_j<0$ for $j\in J_2$ for sufficiently large $k$. 
 Hence, by the optimality condition \eqref{e101} for such $x^k$,  
 $$A_i^T\frac{Ax^k - (\bar b+t_0Aa)}{(\bar\lambda+1/ k)\|Ax^k-(\bar b+t_0Aa)\|}= -1~~\forall i\in   J_1, $$
and
$$A_j^T\frac{Ax^k - (\bar b+t_0Aa)}{(\bar\lambda+1/ k)\|Ax^k-(\bar b+t_0Aa)\|}= 1~~\forall j\in  J_2.$$
Hence, by the definitions of $J_1,J_2$ and $J_{SR}= J_1 \cup J_2$, we have 
 \begin{equation}\label{e612}
    0\neq\frac{Ax^k - (\bar b+t_0Aa)}{(\bar\lambda+1/ k)\|Ax^k-(\bar b+t_0Aa)\|}
    -\frac{A\bar x - \bar b}{\bar\lambda\|A\bar x-\bar b\|}\in \ker(A_{J_{SR}}^T)
 \end{equation}
for all sufficiently large $k$. The inequality holds because these two vectors have different norms.
On the other hand, by the definition of $J_{SR}$, we have 
$$ \left|A_i^T\frac{A(\bar x+t_0a) - (\bar b+t_0Aa)}{\bar\lambda\|A(\bar x+t_0a)-(\bar b+t_0Aa)\|}\right|= \left|A_i^T\frac{A\bar x - \bar b}{\bar\lambda\|A\bar x-\bar b\|}\right|<1~~\forall i\in  J_{SR}^C.$$
Thus $\bar x_{J_{SR}^C}=0$. Furthermore, since $x^k\rightarrow \bar x+t_0a$ as $k\rightarrow \infty$, for sufficiently large $k$, 
$$ \left|A_i^T\frac{Ax^k - (\bar b+t_0Aa)}{(\bar\lambda+1/ k)\|Ax^k-(\bar b+t_0Aa)\|}\right|<1~~\forall i\in J_{SR}^C .$$
By \eqref{e101} for $x^k\in S(\bar b+t_0Aa,\bar\lambda+1/ k)$ implies that $(x^k)_{J_{SR}^C}=0$. Hence, $Ax^k\in \mathrm{rge}(A_{J_{SR}})$. Recalling that $\bar b\in \mathrm{rge}(A_{J_{SR}})$, $a_{J_{SR}^C}=0$ and $\bar x_{J_{SR}^C}=0$, 
it follows that 
 $$0\neq\frac{Ax^k - (\bar b+t_0Aa)}{(\bar\lambda+1/ k)\|Ax^k-(\bar b+t_0Aa)\|}
 -\frac{A\bar x - \bar b}{\bar\lambda\|A\bar x-\bar b\|}\in \mathrm{rge}(A_{J_{SR}}).$$
By \eqref{e612} and the fact $\ker(A_{J_{SR}}^T)\cap \mathrm{rge}(A_{J_{SR}})=\{0\}$, this is impossible. This contradiction implies that $\bar b\notin \mathrm{rge}(A_{J_{SR}})$. The proof is complete.
\end{proof}

\subsection{Single-valued case for SR-LASSO}

This section is devoted to characterizing the locally single-valued and Lipschitz property of solution mapping for the SR-LASSO. To this end, we first answer an open question on solution uniqueness raised by Berk \textit{et al.} \cite{BBH2}. 
Recall that the following condition is proven in \cite{BBH2} to be sufficient for ensuring that $\bar x$ is the unique solution of the SR-LASSO problem \eqref{Pi1}:
\begin{condition}\label{Un2}\cite[Assumption 1]{BBH2}
For $\bar x\in\mathbb{R}^n$ and $I:=\mathrm{supp}(\bar x)$ we have:\\
(i)~$\ker A_I=\{0\}$ and $b\notin \mathrm{rge} A_I;$\\
(ii)~$\exists z\in \ker A_I^T\cap \left\{\frac{b-A\bar x}{\|A\bar x-b\|}\right\}^\perp$ such that $ \left\|A^T_{I^C}\left(\frac{b-A\bar x}{\|A\bar x-b\|}+z \right) \right\|_\infty<\lambda$.
\end{condition}
The open question at the end of section 8 of \cite{BBH2} is: \textit{whether the sufficient Condition \ref{Un2} provided in \cite{BBH2} for the SR-LASSO is also necessary for uniqueness. }

Note that Condition \ref{Un2} requires $\bar x$ to satisfy $b \neq A\bar x$. Therefore, the above question should be restricted to such $\bar x$.
Next, we will show that the answer to the above question is affirmative. 
To do this, we recall the following result, which is contained in the proof of \cite[Theorem 3.4]{BBH2}.

\begin{lemma}(From the proof of \cite[Theorem 3.4]{BBH2})\label{Lu2}
For a solution $\bar x$ of \eqref{Pi1} with $b\neq A\bar x$, let $\bar y=\frac{b-A\bar x}{\|A\bar x-b\|}$ and $\mathcal{S}:=\mathbb{R}_+\{\bar y\}$. Then $\bar x$ is the unique solution of the following problem
\begin{equation}\label{pf}
 \min\limits_{x\in \mathbb{R}^n} \lambda \|x\|_1-\langle A^T\bar y,x\rangle+\delta_{\mathcal{S}}(b-Ax)
\end{equation}
if and only if Condition \ref{Un2} holds.
\end{lemma}

\begin{theorem}\label{Tu3}
If $\bar x$ is the unique solution of the SR-LASSO problem \eqref{Pi1} with $b\neq A\bar x$, then Condition \ref{Un2} holds.
\end{theorem}

\begin{proof}
Let $\bar x \in S(b,\lambda)$. According to Lemma \ref{Lu2}, it suffices to prove that $\bar x$ is the unique solution of \eqref{pf}. Suppose not, then there exists $h\in \mathbb{R}^n\setminus\{0\}$ such that $\bar x+h$ is another solution of \eqref{pf}.  
Since the solution set of \eqref{pf} is convex, we can assume that $\|Ah\|<\|A\bar x-b\|$. Otherwise, we replace $\bar x+h$ by $\bar x+t h$ where $t>0$ is sufficiently small. 
Since both $\bar x$ and $\bar x+h$ are solutions of \eqref{pf}, we have 
\begin{equation}\label{422b}
  b-A(\bar x+h)\in \mathcal{S}=\mathbb{R}_+\{b-A\bar x\}~~\text{and}~~
  \lambda \|\bar x\|_1-\langle A^T\bar y,\bar x\rangle=\lambda \|\bar x+h\|_1-\langle A^T\bar y,\bar x+h\rangle.
\end{equation}
Thus, $Ah\in \mathbb{R}\{b-A\bar x\}$ and $\lambda \|\bar x+h\|_1=\lambda \|\bar x\|_1+\langle A^T\bar y,h\rangle$. Take $t_0\in \mathbb{R}$ such that $Ah= t_0(b-A\bar x)$. Since $\|Ah\|<\|A\bar x-b\|$, we have $|t_0|<1$.
Hence, 
\begin{equation*}
  \aligned
  & \|A(\bar x+h)-b\|+\lambda \|\bar x+h\|_1 \\
  =&(1-t_0)\|A\bar x-b\|+\lambda \|\bar x\|_1+\langle \bar y,Ah\rangle
  \\
  =&(1-t_0)\|A\bar x-b\|+\lambda \|\bar x\|_1+ \left\langle \frac{b-A\bar x}{\|A\bar x-b\|}, t_0(b-A\bar x) \right\rangle
  =\|A\bar x-b\|+\lambda \|\bar x\|_1.
  \endaligned
\end{equation*}
This implies that $\bar x+h \in S(b,\lambda)$, which contradicts the uniqueness of $\bar x$. 
\end{proof}

We are now ready to fully characterize the local single-valued and Lipschitz continuity of solution mapping for the SR-LASSO.

\begin{theorem}\label{ca2a}
Suppose that $\{Ax:x\in S(\bar b,\bar\lambda)\} \neq \{\bar b\}$. 
The solution mapping $S$ of the SR-LASSO
is locally single-valued and Lipschitz at $(\bar b, \bar \lambda)$
if and only if Condition \ref{A2} holds, \textit{i.e.}, $\ker(A_{J_{SR}})=\{0\}$ and $\bar b\notin \mathrm{rge} A_{J_{SR}}$, where $J_{SR}$ is defined as \eqref{J-SR}.
 \end{theorem}
\begin{proof}
It suffices to prove the necessity of Condition \ref{A2}. Since $\{Ax:x\in S(\bar b,\bar\lambda)\} \neq \{\bar b\}$, we can take $\bar x\in S(\bar b,\bar\lambda)$ with $A\bar x\neq b$. 
Denote
$$J_1:=\left\{i\in [n]:A_i^T\frac{A\bar x - \bar b}{\|A\bar x- \bar b\|}=
-\bar \lambda \right\}
~\text{and}~
J_2:= \left\{j\in [n]:A_j^T\frac{A\bar x - \bar b}{\|A\bar x- \bar b\|}= \bar \lambda \right\}, $$
and define $a\in \mathbb R^n$ as \eqref{ea}.
By the assumption that $S$ is locally single-valued and Lipschitz at $(\bar b,\bar \lambda)$, there exists $t>0$ such that $S$ is single-valued at $( \bar b+t Aa,\bar \lambda)$. By Lemma \ref{L613a} (ii), we have $\mathrm{sgn}(\bar x+ta)=\mathrm{sgn}(a)$, $\mathrm{supp}(\bar x+t a)=J_1\cup J_2=J_{SR}$ and $\bar x+t a\in S( \bar b+t Aa, \bar \lambda)$. Since $A(\bar x+t a)\neq \bar b+t Aa$,
by Theorem \ref{Tu3}, we have $\ker(A_{J_{SR}})=\{0\}$,  $\bar b+tAa\notin \mathrm{rge}A_{J_{SR}}$. Recall that $Aa\in \mathrm{rge}A_{J_{SR}}$. We have $\bar b\notin \mathrm{rge}A_{J_{SR}}$.
\end{proof}

\end{document}